\documentclass[11pt,leqno]{amsart}

\usepackage{url}

\usepackage[all,2cell,emtex]{xy}
\usepackage{amssymb,amsmath,bbm,xr,a4wide,color,stmaryrd,enumitem,graphicx}
\usepackage[mathscr]{euscript}
\usepackage{mathrsfs}
\usepackage{textgreek}
\usepackage[colorlinks=true]{hyperref} 

\usepackage[english]{babel}
\DeclareMathOperator{\HGr}{HGr}

\DeclareMathOperator{\sOGr}{\mathbf{GW}}
\DeclareMathOperator{\K}{K}
\DeclareMathOperator{\GW}{GW}
\DeclareMathOperator{\GWe}{GW^{even}_0}

\DeclareMathOperator{\W}{W}
\DeclareMathOperator{\Sym}{Sym}
\DeclareMathOperator{\id}{id}

\DeclareMathOperator{\im}{im}
\DeclareMathOperator{\rank}{rank}
\DeclareMathOperator{\Hom}{Hom}
\DeclareMathOperator{\Ext}{Ext}
\DeclareMathOperator{\Spec}{Spec}
\DeclareMathOperator{\fl}{fl}
\DeclareMathOperator{\modul}{mod}
\DeclareMathOperator{\VB}{VB}
\DeclareMathOperator{\Db}{D^b}
\DeclareMathOperator{\Dbfl}{D^b_{fl}}
\DeclareMathOperator{\Sp}{Sp}
\DeclareMathOperator{\BO}{\mathbf{BO}}

\newcommand{\NN} {\mathbb N}
\newcommand{\ZZ} {\mathbb Z}
\newcommand{\QQ} {\mathbb Q}
\newcommand{\Zh} {\ZZ[\frac{1}{2}]}
\newcommand{\Sy} {\mathfrak{S}}
\newcommand{\OO} {\mathcal O}
\renewcommand{\AA} {\mathbb A}
\newcommand{\PP} {\mathbb P}

\newcommand{\HP} {\mathrm{H}\mathbb{P}}
\newcommand{\Tc} {\mathcal{T}}
\newcommand{\Ec} {\mathcal{E}}
\newcommand{\Gc} {\mathcal{G}}
\newcommand{\SH} {\mathcal{SH}}
\newcommand{\Ho} {\mathcal{H}}
\newcommand{\GWi}[2]{\GW^{#1}\{\frac{1}{{#2}^*}\}}
\newcommand{\GWibig}[2]{\GW^{#1}\Big\{\frac{1}{{#2}^*}\Big\}}

\newcommand{\ext}{\mathchoice{{\textstyle\bigwedge}}
    {{\bigwedge}}
    {{\textstyle\wedge}}
    {{\scriptstyle\wedge}}}

\makeatletter
\newcommand{\bigperp}{%
  \mathop{\mathpalette\bigp@rp\relax}%
  \displaylimits
}

\newcommand{\bigp@rp}[2]{%
  \vcenter{
    \m@th\hbox{\scalebox{\ifx#1\displaystyle2.1\else1.5\fi}{$#1\perp$}}
  }%
}

\title{The stable Adams operations on Hermitian $K$-theory}

\author{Jean Fasel}
\address{Institut Fourier - UMR 5582, Universit\'e Grenoble-Alpes, CS 40700, 38058 Grenoble Cedex 9, France}
\email{jean.fasel@univ-grenoble-alpes.fr}
\urladdr{\url{https://www-fourier.univ-grenoble-alpes.fr/~faselj/}}

\author{Olivier Haution}
\address{Mathematisches Institut, Ludwig-Maximilians-Universit\"at M\"unchen, Theresienstr.\ 39, D-80333 M\"unchen, Germany}
\email{olivier.haution@gmail.com}
\urladdr{\url{https://haution.gitlab.io/}}
\thanks{This work was supported by the DFG grant HA 7702/5-1 and Heisenberg fellowship HA 7702/4-1.}

\date{\today}

\newtheorem{thm}{Theorem}[subsection]
\newtheorem{prop}[thm]{Proposition}
\newtheorem{lem}[thm]{Lemma}
\newtheorem{cor}[thm]{Corollary}
\newtheorem{prop-df}[thm]{Proposition--Definition}

\theoremstyle{remark} 
\newtheorem{rem}[thm]{Remark}

\theoremstyle{definition} 
\newtheorem{df}[thm]{Definition}
\newtheorem{notation}[thm]{Notation}

\renewcommand{\theequation}{\thesubsection.\alph{equation}}

\begin{document}

\begin{abstract}
We prove that exterior powers of (skew-)symmetric bundles induce a $\lambda$-ring structure on the ring $\GW^0(X) \oplus \GW^2(X)$, when $X$ is a scheme where $2$ is invertible. Using this structure, we define stable Adams operations on Hermitian $K$-theory. As a byproduct of our methods, we also compute the ternary laws associated to Hermitian $K$-theory.
\end{abstract}

\maketitle

\setcounter{tocdepth}{2}
\tableofcontents

\section*{Introduction}

From their introduction by Adams in his study of vector fields on spheres \cite{Adams62}, Adams operations have been extremely useful in solving various problems in topology, algebra and beyond. One may mention for instance the proof of Serre vanishing conjecture by Gillet-Soul\'e \cite{Gillet87}, or their use in intersection theory. In algebraic geometry, the work of several authors permitted to extend these operations (initially defined at the level of the Grothendieck group $K_0$) to the whole world of $K$-theory; the most recent and probably most natural extension being due to Riou \cite{Riou10} using (stable) motivic homotopy theory. 

Over a scheme $X$, it is often useful to study vector bundles endowed with some extra decoration, such as a symmetric or a symplectic form. The analogues of the Grothendieck group $K_0(X)$ in this context are the so-called Grothendieck--Witt groups (or Hermitian $K$-theory groups) $\GW^i(X)$ for $i\in \ZZ/4$ (see e.g.\ \cite{Schlichting17}), which classify symmetric and symplectic bundles \cite{Walter03}. Very often, the constructions and questions pertaining to algebraic $K$-theory can be generalized to the context of Grothendieck--Witt groups. For instance, Serre's Vanishing Conjecture makes sense in this broader context \cite{Fasel08c}.\\

As for the Adams operations, Zibrowius \cite{Zibrowius15,Zibrowius18} has proved that the exterior power operations on symmetric bundles yield a $\lambda$-ring structure on the Grothendieck--Witt group $\GW^0(X)$ of any smooth variety $X$ over a field of characteristic not two. This provides in particular Adams operations on these groups. It is not very difficult to construct $\lambda$-operations in $\GW^0(X)$, and a significant portion of the papers \cite{Zibrowius15,Zibrowius18} consists in showing that this pre-$\lambda$-ring is actually a $\lambda$-ring, which means that the $\lambda$-operations verify certain additional relations pertaining to their multiplicative and iterative behaviour. In particular, it is not so difficult to construct the Adams operations $\psi^n$, but much harder to show that they are multiplicative and verify the relations $\psi^{mn}=\psi^m \circ \psi^n$. To prove that $\GW^0(X)$ is a $\lambda$-ring, Zibrowius followed the strategy used in \cite{SGA6} for the analog problem in $K$-theory, and reduced the question to proving that the symmetric representation ring $\GW^0(G)$ of an affine algebraic group $G$ (over a field of characteristic not two) is a $\lambda$-ring. This is done by further reducing to the case when $G$ is the split orthogonal group, and using explicit descriptions of the representations of certain subgroups in that case.

A first purpose of this paper is to extend the construction of Zibrowius in two directions:
\begin{enumerate}[label={(\arabic*)}]
\item allow $X$ to be an arbitrary quasi-compact quasi-separated $\Zh$-scheme admitting an ample family of line bundles,

\item replace $\GW^0(X)$ with $\GW^{\pm}(X)$, the ring of symmetric and symplectic forms.
\end{enumerate}
The objective is achieved by first showing that the map $\GW^0(G) \to \GW^0(G_{\QQ})$ is injective, when $G$ is a split reductive algebraic group over $\Zh$. Since the target is a $\lambda$-ring by the result of Zibrowius, so is $\GW^0(G)$, and thus also $\GW^0(X)$ when $X$ is as in (1).

For (2), a natural strategy is to mimic Zibrowius's proof, by considering not just symmetric representations of algebraic groups, but also skew-symmetric ones. Although we believe that this idea might work, we were not able to implement it satisfyingly. Instead we observe that we may pass from $\GW^-(X)$ to $\GW^+(X)$ using the quaternionic projective bundle theorem \cite{Panin10pred}.\\

The Witt groups are natural companions of the Grothendieck--Witt groups, obtained from them by modding out the hyperbolic classes. Their behaviour is somewhat easier to understand, and they keep track of an important part of the quadratic information, while forgetting some of the $K$-theoretic information. Our $\lambda$-ring structure on the Grothendieck--Witt groups does not descend to one on the Witt groups. There is a good reason for this: the Witt ring cannot admit a (functorial) $\lambda$-ring structure, because it takes the value $\mathbb{F}_2$ on every algebraically closed field, and $\mathbb{F}_2$ has no such structure. Nonetheless, we prove that the odd Adams operations (as well as the even ones when additionally $-1$ is a square) do descend to operations on the Witt ring. It would be interesting to find algebraic axioms describing a weak form of the structure of $\lambda$-ring (including odd Adams operations) which applies to the Witt ring, but we will not investigate this question further in this paper.\\

The next natural step consists in considering the groups $\GW^i(X)$ for $i$ odd, as well as the higher Grothendieck--Witt groups $\GW^i_j(X)$ for $j\in \ZZ$. To do so, we focus on Adams operations, and follow the approach pioneered by Riou \cite{Riou10} to construct stable versions of those. The fact that $\GW^{\pm}(X)$ is a $\lambda$-ring ends up being a crucial input, allowing us to understand the behaviour of the Adams operations with respect to stabilization. This approach is carried out in Section \ref{sec:stable}, where we build a morphism of motivic ring spectra, for any integer $n\in \NN$
\[
\Psi^n\colon\sOGr\to \sOGr\Big[\frac{1}{n^*}\Big].
\]
Here the left-hand side is the spectrum representing Hermitian $K$-theory and the right-hand side is the same after inversion of the class $n^*\in\GW^+(X)$, which equals $n$ when $n$ is odd, and the class of the hyperbolic $n$-dimensional symmetric form when $n$ is even. These operations extend the Adams operations on $K$-theory, in the sense that there is a commutative diagram of motivic ring spectra
\[
\xymatrix{\sOGr\ar[r]^-{\Psi^n}\ar[d] & \sOGr[\frac{1}{n^*}] \ar[d]\\
\mathbf{BGL}\ar[r]^-{\Psi^n} & \mathbf{BGL}[\frac{1}{n}]}
\]
in which the vertical morphisms are the forgetful maps and the bottom horizontal morphism is the Adams operation on $K$-theory defined by Riou \cite[Definition 5.3.2]{Riou10}.

When $n$ is even, inverting $n^*$ in $\GW^+(X)$ seems to be a fairly destructive procedure, so in practice the stable even Adams operations are unlikely to be very valuable improvements of their $K$-theoretic counterparts. By contrast, we expect that the odd operations will be useful in many situations. For instance, Bachmann and Hopkins recently used them in \cite{Bachmann20a} to compute the $\eta$-inverted homotopy sheaves of the algebraic symplectic and special linear cobordism spaces. Their construction of Adams operations is quite different in spirit to the one presented here but satisfy (almost) the same properties (see \cite[Remark 3.2]{Bachmann20a}).

In the last section of this paper, we offer an application under the form of the computation of the ternary laws associated to Hermitian $K$-theory. These laws are the analogue,  in the context of Sp-oriented ring spectra, of the formal group laws associated to any oriented ring spectrum. In short, they express the characteristic classes of a threefold product of symplectic bundles of rank $2$, and are expected to play an important role in the classification of Sp-oriented cohomology theories. We refer the interested reader to \cite{Deglise20} for more information on these laws.\\

\noindent \textbf{Acknowledgments}.
The first named author is grateful to Aravind Asok, Baptiste Calm\`es and Fr\'ed\'eric D\'eglise for useful discussions. Both authors warmly thank Alexey Ananyevskiy for sharing a preprint on Adams operations which has been a source of inspiration for the results of the present paper, and Tom Bachmann for very useful suggestions. They also heartily thank the referee for a careful reading and useful comments that helped correct mistakes and improve the exposition.

\section{Grothendieck--Witt groups and spectra}\label{sec:GWspectra}
\label{sec:spectrum}
\numberwithin{thm}{section}
\numberwithin{equation}{section}
\renewcommand{\theequation}{\thesection.\alph{equation}}
All schemes will be assumed to be quasi-compact and quasi-separated, and to admit an ample family of line bundles.

Let $X$ be a scheme. In this paper, we will denote by $\GW^+(X)$, resp.\ $\GW^-(X)$, the Grothendieck--Witt group of symmetric forms, resp.\ skew-symmetric forms, defined e.g.\ in \cite[\S6]{Walter03} using the exact category of vector bundles over $X$. The product of two skew-symmetric forms being symmetric, we have a pairing
\[
\GW^-(X)\times \GW^-(X)\to \GW^+(X)
\]
turning $\GW^{\pm}(X)=\GW^+(X)\oplus  \GW^-(X)$ into a (commutative) $\ZZ/2$-graded ring. 

Assume now that $X$ is a scheme over $\Zh$. Following \cite[Definition 9.1]{Schlichting17}, we can consider the Grothendieck--Witt groups $\GW^i_j(X)$ for any $i,j\in \ZZ$ which are $4$-periodic in $i$ in the sense that there are natural isomorphisms $\GW^i_j(X)\simeq \GW^{i+4}_j(X)$ for any $i\in\ZZ$. For $X$ affine and $i=0$, the groups $\GW_j^0(X)$ are (naturally isomorphic to) the orthogonal $K$-theory groups $\mathrm{KO}_j(X)$ as defined by Karoubi, while for $i=2$ (and $X$ still affine) the groups $\GW_j^2(X)$ are (naturally isomorphic to) the symplectic $K$-theory groups $\mathrm{KSp}_j(X)$ (\cite[Corollary A.2]{Schlichting17}). Also by \cite[Theorem~6.1]{Walter03} and \cite[Proposition~5.6]{Schlichting17} we have natural isomorphisms $\GW^+(X) \simeq \GW^0_0(X)$ and $\GW^-(X) \simeq \GW^2_0(X)$. 

\begin{notation}
\label{not:tau}
We will denote by $h \in \GW^0_0(\Spec(\Zh))$, resp.\ $\tau\in \GW^2_0(\Spec(\Zh))$, the class of the hyperbolic symmetric, resp.\ skew-symmetric, bilinear form. When $u \in (\Zh)^\times$, we will denote by $\langle u \rangle \in \GW^0_0(\Spec(\Zh))$ the class of the symmetric bilinear form $(x,y) \mapsto uxy$, and write $\epsilon = -\langle -1 \rangle$. Thus $h=1 - \epsilon$.
\end{notation}

The collection of groups $\GW^i_j(X)$ fit into a well-behaved cohomology theory, which is $\mathrm{SL}^c$-oriented by \cite[Theorem~5.1]{Panin10prec}, and in particular $\Sp$-oriented. The functors $X\mapsto \GW^i_j(X)$ are actually representable by explicit (geometric) spaces $\GW^i$ in the $\AA^1$-homotopy category $\Ho(\Zh)$ of Morel-Voevodsky (see \cite[Theorem 1.3]{Schlichting13})
\[
[\Sigma^j_{S^1}X_+,\GW^i]_{\Ho(\Zh)}=\GW^i_j(X).
\]
Further, one can express the aforementioned periodicity under the following form: there exists an element $\gamma\in \GW^4_0(\Spec(\Zh))$ such that multiplication by $\gamma$ induces the periodicity isomorphisms
\begin{equation}
\label{eq:periodicity}
\GW^i\simeq \GW^{i+4}.
\end{equation}

When $X$ is a $\Zh$-scheme, the $\ZZ$-graded ring
\begin{equation}
\label{eq:GWe}
\GWe(X):= \bigoplus_{j\in\ZZ}\GW_0^{2j}(X)
\end{equation}
can be identified with the $\ZZ$-graded subring $\widehat{\GW^{\pm}(X)}$ of $\GW^{\pm}(X)[x^{\pm 1}]$ defined in Appendix~\ref{sec:gradedrings} (where $\gamma$ corresponds to $x^2$), and we have a canonical isomorphism of $\ZZ/2$-graded rings $\GWe(X)/(\gamma -1) \simeq \GW^{\pm}(X)$.

The $\PP^1$-projective bundle theorem of Schlichting \cite[Theorem 9.10]{Schlichting17} allows to build a ring spectrum $\sOGr$ in $\SH(\Zh)$, having the property to represent Hermitian $K$-theory. A convenient construction is recalled in \cite[Theorem 12.2]{Panin10prec}, and we explain the relevant facts in the next few lines in order to fix notations. 

Recall first that Panin and Walter \cite{Panin10pred} defined a smooth affine $\Zh$-scheme $\HP^n$ for any $n\in\NN$, called the quaternionic projective space. On $\HP^n$, there is a canonical bundle $U$ of rank $2$ endowed with a symplectic form $\varphi$, yielding a canonical element $u=(U,\varphi)\in \GW^-(\HP^n)$. For any $n\in\NN$, there are morphisms
\begin{equation}
\label{eq:i_n}
i_n\colon\HP^n\to \HP^{n+1}
\end{equation}
such that $i_n^*u=u$, whose colimit (say in the category of sheaves of sets) is denoted by $\HP^{\infty}$. It is a geometric model of the classifying space $\mathrm{BSp}_2$ of rank $2$ symplectic bundles. As $\HP^0=\Spec(\Zh)$, we consider all these schemes as pointed by $i_0$ and note that $i_0^*(u)=\tau$. Recall moreover from \cite[Theorem~9.8]{Panin10prec} that $\HP^1$ is $\mathbb{A}^1$-weak equivalent to $(\PP^1)^{\wedge 2}$. In fact $\HP^1=Q_4$, where the latter is the affine scheme considered for instance in \cite{Asok14}. 

\begin{notation}
We set $\Tc:=\HP^1$, that we consider as pointed by $i_0$. We also denote by $\Omega_\Tc$ the right adjoint of the endofunctor $\Tc\wedge(-)$ of $\Ho(\Zh)$.
\end{notation}

The spectrum $\sOGr$ is defined as the $\Tc$-spectrum whose component in degree $n$ is $\GW^{2n}$ and bonding maps
\begin{equation}
\label{eq:def_sigma}
\sigma \colon \Tc\wedge \GW^{2n}\to \GW^{2n+2}
\end{equation}
induced by multiplication by the class $u-\tau$ in $\GW^2_0(\Tc)$. This $\Tc$-spectrum determines uniquely a $\PP^1$-spectrum in view of \cite[Proposition 2.22]{Riou07} or \cite[Theorem 12.1]{Panin10prec}, which has the property that 
\[
\GW^{i}_j(X)=[\Sigma_{\PP^1}^{\infty}X_+,\Sigma_{S^1}^{-j}\Sigma_{\PP^1}^i\sOGr]_{\SH(\Zh)}
\]
for a smooth $\Zh$-scheme $X$.

If now $X$ is a regular $\Zh$-scheme with structural morphism $p_X\colon X\to \Spec(\Zh)$, we can consider the functor $p_X^*\colon \SH(\Zh)\to \SH(X)$ and the spectrum $p_X^*\sOGr$. On the other hand, one can consider the $\PP^1_X$-spectrum $\sOGr_X$ representing Grothendieck--Witt groups in the stable category $\SH(X)$. It follows from \cite[discussion before Theorem 13.5]{Panin10prec} that the natural map $p_X^*\sOGr\to \sOGr_X$ is in fact an isomorphism. Consequently,
\[
\GW^i_j(X)=[\Sigma_{\PP^1}^{\infty}X_+,\Sigma_{S^1}^{-j}\Sigma_{\PP^1}^ip_X^*\sOGr]_{\SH(X)}
\]
and we say that $\sOGr$ is an \emph{absolute $\PP^1$-spectrum} over $\Zh$. It is in fact an absolute \emph{ring spectrum} by \cite[Theorem 13.4]{Panin10prec}.

\section{Exterior powers and rank two symplectic bundles}

When $V$ is a vector bundle on a scheme $X$, we denote its dual by $V^\vee$. A bilinear form on $V$ is a morphism of vector bundles $\nu \colon V \to V^\vee$. When $x,y \in H^0(X,V)$, we will sometimes write $\nu(x,y)$ instead of $\nu(x)(y)$. We will abuse notation, and denote by $\ext^n\nu$, for $n\in \NN$, the bilinear form on $\ext^nV$ given by the composite $\ext^nV \xrightarrow{\ext^n\nu} \ext^n(V^\vee) \to (\ext^nV)^\vee$. We will also denote the pair $(\ext^nV,\ext^n\nu)$ by $\ext^n(V,\nu)$. Similar conventions will be used for the symmetric or tensor powers of bilinear forms, or their tensor products. 

Explicit formulas for symmetric and exterior powers are given as follows. Let $n$ be an integer, and denote by $\Sy_n$ the symmetric group on $n$ letters and by $\epsilon \colon \Sy_n \to \{-1,1\}$ the signature homomorphism. Then for any open subscheme $U$ of $X$ and $x_1,\dots,x_n,y_1,\dots,y_n \in H^0(U,V)$, we have
\begin{equation}
\label{eq:explicit_sym}
(\Sym^n\nu)(x_1\cdots x_n,y_1 \cdots y_n)= \sum_{\sigma \in \Sy_n} \nu(x_1,y_{\sigma(1)})\cdots\nu(x_n,y_{\sigma(n)}),
\end{equation}
\begin{equation}
\label{eq:explicit_ext}
(\ext^n\nu)(x_1\wedge \cdots \wedge x_n,y_1\wedge \cdots \wedge y_n)= \sum_{\sigma \in \Sy_n} \epsilon(\sigma) \nu(x_1,y_{\sigma(1)})\cdots\nu(x_n,y_{\sigma(n)}),
\end{equation}
or more succinctly
\begin{equation}
\label{eq:explicit_ext:det}
(\ext^n\nu)(x_1\wedge \cdots \wedge x_n,y_1\wedge \cdots \wedge y_n)= \det(\nu(x_i,y_j)).
\end{equation}

If $V,W$ are vector bundles equipped with bilinear forms $\nu,\mu$, then for any $i,j$ the bilinear form $\ext^{i+j}(\nu \perp \mu)$ restricts to $(\ext^i\nu) \otimes (\ext^j\mu)$ on $(\ext^i V) \oplus (\ext^j W) \subset \ext^{i+j}(V \oplus W)$. This yields an isometry, for any $n \in \NN$
\begin{equation}
\label{eq:ext_sum}
\ext^n(V \oplus W,\nu \oplus \mu) \simeq \bigperp_{i=0}^n \ext^i(V,\nu) \otimes \ext^{n-i}(W,\mu)
\end{equation}

\begin{lem}
\label{lemm:det_tensor}
Let $(E,\varepsilon)$ and $(F,\varphi)$ be vector bundles over a scheme $X$ equipped with bilinear forms, of respective ranks $e$ and $f$. Then we have an isometry
\[
(\ext^e(E,\varepsilon))^{\otimes f} \otimes (\ext^f(F,\varphi))^{\otimes e} \simeq \ext^{ef}(E \otimes F,\varepsilon \otimes \varphi).
\]
\end{lem}
\begin{proof}
Let us first assume that $E,F$ are free and that $X=\Spec R$ is affine. Let $(x_1,\dots,x_e)$, resp.\ $(y_1,\dots,y_f)$, be an $R$-basis of $H^0(X,E)$, resp.\ $H^0(X,F)$. Then the element
\begin{equation}
\label{eq:def_z}
z=(x_1 \wedge \dots \wedge x_e)^{\otimes f} \otimes (y_1 \wedge \dots \wedge y_f)^{\otimes e}
\end{equation}
is a basis of $H^0(X,(\ext^e E)^{\otimes f} \otimes (\ext^f F)^{\otimes e})$, and the element
\begin{equation}
\label{eq:def_u}
u =  (x_1 \otimes y_1) \wedge \dots \wedge (x_1 \otimes y_f) \wedge (x_2 \otimes y_1) \wedge \dots \wedge (x_e \otimes y_f)
\end{equation}
is a basis of $H^0(X,\ext^{ef} (E\otimes F)$. The mapping  $z \mapsto u$ then defines an isomorphism of line bundles
\begin{equation}
\label{eq:loc_isom}
(\ext^e E)^{\otimes f} \otimes (\ext^f F)^{\otimes e} \xrightarrow{\sim} \ext^{ef} (E\otimes F),
\end{equation}
Consider now the matrices
\[
A=(\varepsilon(x_i,x_j)) \in M_e(R), \quad B=(\varphi(y_i,y_j)) \in M_f(R).
\]
By \eqref{eq:explicit_ext:det} we have
\[
((\ext^e \varepsilon)^{\otimes f} \otimes (\ext^f F)^{\otimes e})(z,z)=(\det A)^f \cdot (\det B)^e,
\]
and $\ext^{ef} (\varepsilon \otimes \varphi)(u,u)$ is the determinant of the block matrix
\[
C=\left(\begin{matrix}
\varepsilon(x_1,x_1) B & \dots & \varepsilon(x_1,x_e)B\\
\vdots &&\vdots\\
\varepsilon(x_e,x_1) B & \dots & \varepsilon(x_e,x_e)B
\end{matrix}\right) \in M_{ef}(R).
\]
It then follows from \cite[III, \S9, Lemme 1, p.112]{Bou-A-13} that
\[
\det C = \det(\det(\varepsilon(x_i,x_j) B)) = \det(\det(A)B^e)=(\det A)^f \cdot (\det B)^e.
\]
Therefore
\[
((\ext^e \varepsilon)^{\otimes f} \otimes (\ext^f F)^{\otimes e})(z,z) = \ext^{ef} (\varepsilon \otimes \varphi)(u,u),
\]
which shows that \eqref{eq:loc_isom} is the required isometry.

Next, assume given $R$-linear automorphisms $\alpha \colon H^0(X,E) \to H^0(X,E)$ and $\beta \colon H^0(X,F) \to H^0(X,F)$. Replacing the basis $(x_1,\dots,x_e)$ and $(y_1,\dots,y_f)$ by their images under $\alpha$ and $\beta$ multiplies the element \eqref{eq:def_z} by the quantity $(\det \alpha)^e \cdot (\det \beta)^f$, and the element \eqref{eq:def_u} by the same quantity (this is a similar determinant computation as above, based on \cite[III, \S9, Lemme 1, p.112]{Bou-A-13}). We deduce that the isometry \eqref{eq:loc_isom} glues when $E,F$ are only (locally free) vector bundles, and $X$ is possibly non-affine.
\end{proof}

\begin{lem}
\label{lem:ext_n-1}
Let $V$ be a vector bundle of constant rank $n$ over a scheme $X$, equipped with a nondegenerate bilinear form $\nu$. Then we have an isometry
\[
\ext^{n-1} (V,\nu) \simeq (V,\nu) \otimes \ext^n(V,\nu).
\]
\end{lem}
\begin{proof}
The natural morphism $(\ext^{n-1}V) \otimes V \to \ext^n V$ induces a morphism $\ext^{n-1} V \to \Hom(V,\ext^n V)$. As $V$ is a vector bundle (of finite rank) the natural morphism $V^\vee \otimes \ext^n V \to \Hom(V,\ext^n V)$ is an isomorphism. Composing with the inverse of $\nu \otimes \id_{\ext^n V}$, we obtain a morphism
\[
s \colon \ext^{n-1} V \to V \otimes \ext^n V.
\]
To verify that it induces the required isometry, we may argue locally and assume that $V$ is free and $X =\Spec R$ is affine. Pick an $R$-basis $(v_1,\dots,v_n)$ of $H^0(X,V)$. Then $(w_1,\dots,w_n)$ is an $R$-basis of $H^0(X,\ext^{n-1}V)$, where $w_i = (-1)^{n-i} v_1 \wedge \dots \wedge \widehat{v_i} \wedge v_n$. Let $z= v_1 \wedge \dots \wedge v_n \in H^0(X,\ext^nV)$, and note that $w_i \wedge v_i = z$ for all $i\in \{1,\dots,n\}$. Consider the unique elements $v_1^*,\dots,v_n^* \in H^0(X,V)$ satisfying $\nu(v_i^*,v_j) = \delta_{ij}$ (Kronecker symbol) for all $i,j \in \{1,\dots,n\}$. Then we have
\begin{equation}
\label{eq:s_w_v_z}
s(w_i) = v_i^* \otimes z, \quad \text{ for $i=1,\dots,n$}.
\end{equation}
Consider the matrix $A=(\nu(v_i,v_j)) \in M_n(R)$. Observe that the $j$-th coordinate of $v_i^*$ in the basis $(v_1,\dots,v_n)$ is the $(i,j)$-th coefficient of the matrix $A^{-1}$, from which it follows that
\begin{equation}
\label{eq:t_A_-1}
\!^{t}(A^{-1}) = (\nu(v_i^*,v_j^*)) \in M_n(R).
\end{equation}
Let $k,l \in \{1,\dots,n\}$. It follows from \eqref{eq:explicit_ext:det} that $(\ext^{n-1} \nu)(w_k,w_l)$ is the $(k,l)$-th cofactor of the matrix $A$, and thus coincides with the $(k,l)$-th coefficient of the matrix $(\det A) \cdot \!^{t}(A^{-1})$. In view of \eqref{eq:t_A_-1}, we deduce that (using \eqref{eq:explicit_ext:det} for the last equality)
\[
(\ext^{n-1} \nu)(w_k,w_l) = \nu(v_k^*,v_l^*) \cdot \det A = \nu(v_k^*,v_l^*) \cdot (\ext^n\nu)(z,z).
\]
By the formula \eqref{eq:s_w_v_z}, this proves that $s$ is the required isometry.
\end{proof}

In the rest of the section, we fix a $\Zh$-scheme $X$. By a symplectic bundle on $X$, we will mean a vector bundle on $X$ equipped with a nondegenerate skew-symmetric form. For an invertible element $\lambda \in H^0(X,\OO_X)$, we denote by $\langle \lambda \rangle$ the trivial line bundle on $X$ equipped with the nondegenerate bilinear form given by $(x,y) \mapsto \lambda xy$.

\begin{lem}
\label{lemm:det_symplectic}
Let $(V,\nu)$ be a symplectic bundle of constant rank $n$ over $X$. Then the exists an isometry
\[
\ext^n(V,\nu) \simeq \langle 1 \rangle.
\]
\end{lem}
\begin{proof}
We may assume that $X \neq \varnothing$ and $n \geq 1$. Then we may write $n=2m$ for some integer $m$ (the form induced by $(V,\nu)$ over the residue field of a closed point of $X$ is skew-symmetric, hence symplectic as $2$ is invertible, and such forms over fields have even dimension \cite[I, (3.5)]{Milnor73}). The morphism 
\[
V^{\otimes n} \simeq V^{\otimes m} \otimes V^{\otimes m} \to  \ext^mV \otimes \ext^mV \xrightarrow{\ext^m\nu \otimes \id} (\ext^mV)^\vee \otimes \ext^mV \to \OO_X
\]
descends to a morphism $\lambda_{(V,\nu)} \colon \ext^nV \to \OO_X$. If $(V_i,\nu_i)$, for $i=1,2$, are symplectic bundles over $X$ of ranks $n_i=2m_i$ such that $(V,\nu) = (V_1,\nu_1) \perp (V_2,\nu_2)$, we have a commutative diagram
\[
\xymatrix{
\ext^m V\ar[rrr]^{\ext^m\nu} &&& (\ext^m V)^\vee\\ 
(\ext^{m_1}V_1) \otimes (\ext^{m_2}V_2) \ar[rrr]^{\ext^{m_1}\nu_1 \otimes \ext^{m_2}\nu_2}\ar[u] &&& (\ext^{m_1}V_1)^\vee \otimes (\ext^{m_2}V_2)^\vee \ar[u]
}
\]
Therefore the identification $\ext^nV = \ext^{n_1}V_1 \otimes \ext^{n_2} V_2$ yields an identification $\lambda_{(V,\nu)} = \lambda_{(V_1,\nu_1)} \otimes \lambda_{(V_2,\nu_2)}$.

In order to prove that $\lambda_{(V,\nu)}$ induces the claimed isometry, we may assume that $X$ is the spectrum of a local ring. In this case the nondegenerate skew-symmetric form $(V,\nu)$ is hyperbolic \cite[I, (3.5)]{Milnor73}. Given the behaviour of $\lambda_{(V,\nu)}$ with respect to orthogonal sums, we may assume that $n=2$ and that $(V,\nu)$ is the hyperbolic plane. So there exists a basis $(v_1,v_2)$ of $H^0(X,V)$ such that 
\[
\nu(v_1,v_1)=0, \; \nu(v_2,v_2)=0 \text{ and } \nu(v_1,v_2)=1.
\]
By \eqref{eq:explicit_ext} we have
\[
(\ext^2\nu)(v_1 \wedge v_2,v_1 \wedge v_2) = 1.
\]
Since $\lambda_{(V,\nu)}(v_1 \wedge v_2) = \nu(v_1,v_2)=1 \in H^0(X,\OO_X)$, it follows that $\lambda_{(V,\nu)}$ induces an isometry $\ext^2(V,\nu) \simeq \langle 1 \rangle$.
\end{proof}

Let $V$ be a vector bundle over $X$. Consider the involution $\sigma$ of $V^{\otimes 2}$ exchanging the two factors. Set $V^{\otimes 2}_+ = \ker(\sigma -\id)$ and $V^{\otimes 2}_- = \ker(\sigma +\id)$. Since $2$ is invertible we have a direct sum decomposition $V^{\otimes 2} =V^{\otimes 2}_+ \oplus V^{\otimes 2}_-$.

Let now $\nu$ be a bilinear form on $V$. There are induced bilinear forms  $\nu^{\otimes 2}_+$ on $V^{\otimes 2}_+$ and $\nu^{\otimes 2}_-$ on $V^{\otimes 2}_-$. Writing $(V,\nu)^{\otimes 2}_+$, resp.\ $(V,\nu)^{\otimes 2}_-$, instead of $(V^{\otimes 2}_+,\nu^{\otimes 2}_+)$, resp.\ $(V^{\otimes 2}_-,\nu^{\otimes 2}_-)$, we have an orthogonal decomposition
\begin{equation}
\label{eq:decomp_Votimes2}
(V,\nu)^{\otimes 2} = (V,\nu)^{\otimes 2}_+ \perp   (V,\nu)^{\otimes 2}_-.
\end{equation}

\begin{lem}
\label{lem:+-_SymLambda}
There are isometries
\[
(V,\nu)^{\otimes 2}_+ \simeq \langle 2 \rangle \otimes \Sym^2(V,\nu) \quad \text{and} \quad (V,\nu)^{\otimes 2}_- \simeq \langle 2 \rangle \otimes \ext^2(V,\nu).
\]
\end{lem}
\begin{proof}
It is easy to see that the morphism
\[
i \colon \ext^2V \to  V^{\otimes 2} \quad ; \quad v_1 \wedge v_2 \mapsto v_1 \otimes v_2 - v_2 \otimes v_1,
\]
induces an isomorphism $\ext^2V \simeq V^{\otimes 2}_-$. If $U$ is an open subscheme of $X$ and $v_1,v_2,w_1,w_2 \in H^0(U,V)$, we have, using \eqref{eq:explicit_ext}
\begin{align*}
&\nu^{\otimes 2}(i(v_1\wedge v_2), i(w_1\wedge w_2))\\
=&\nu^{\otimes 2}(v_1 \otimes v_2-v_2\otimes v_1,w_1\otimes w_2 - w_2 \otimes w_1) \\ 
=& \nu(v_1,w_1)\nu(v_2,w_2) - \nu(v_2,w_1)\nu(v_1,w_2) - \nu(v_1,w_2)\nu(v_2,w_1) +\nu(v_2,w_2)\nu(v_1,w_1)\\
=&2\nu(v_1,w_1)\nu(v_2,w_2) - 2\nu(v_2,w_1)\nu(v_1,w_2)\\
=&2 (\ext^2\nu)(v_1\wedge v_2,w_1\wedge w_2),
\end{align*}
proving the second statement. The first is proved in a similar fashion, using the morphism
\[
\Sym^2V \to  V^{\otimes 2} \quad ; \quad v_1 v_2 \mapsto v_1 \otimes v_2 + v_2 \otimes v_1.\qedhere
\]
\end{proof}

\begin{lem}
\label{lem:tens2_decomp}
There is an isometry
\[
(V,\nu)^{\otimes 2} \simeq \langle 2 \rangle \otimes \Big(\Sym^2(V,\nu) \perp \ext^2(V,\nu) \Big).
\]
\end{lem}
\begin{proof}
This follows from Lemma \ref{lem:+-_SymLambda} and \eqref{eq:decomp_Votimes2}.
\end{proof}

\begin{lem}
\label{lem:lambda_EF}
Let $E,F$ be vector bundles over $X$, respectively equipped with bilinear forms $\varepsilon,\varphi$. Then there is an isometry
\[
\ext^2(E \otimes F,\varepsilon \otimes \varphi) \simeq  \Big(\langle 2 \rangle \otimes \Sym^2(E,\varepsilon) \otimes \ext^2(F,\varphi)\Big) \perp  \Big(\langle 2\rangle \otimes \ext^2(E,\varepsilon) \otimes \Sym^2(F,\varphi)\Big).
\]
\end{lem}
\begin{proof}
It is easy to see that there is an isometry
\[
(E\otimes F,\varepsilon \otimes \varphi)^{\otimes 2}_- \simeq \Big((E,\varepsilon)^{\otimes 2}_+ \otimes (F,\varphi)^{\otimes 2}_-\Big) \perp \Big((E,\varepsilon)^{\otimes 2}_- \otimes (F,\varphi)^{\otimes 2}_+\Big),
\]
so that the statement follows by five applications of Lemma \ref{lem:+-_SymLambda} (and tensoring by the form $\langle 2^{-1} \rangle$).
\end{proof}

\begin{prop}
\label{prop:lambda_22}
Let $E,F$ be rank two vector bundles over a $\Zh$-scheme $X$, equipped with nondegenerate skew-symmetric forms $\varepsilon, \varphi$. Then we have in $\GW^+(X)$:
\[
[\ext^n(E \otimes F,\varepsilon \otimes \varphi)]=
\begin{cases}
[(E,\varepsilon) \otimes (F,\varphi)] &\text{if $n \in \{1,3\}$,}\\
[(E,\varepsilon)^{\otimes 2}] + [(F,\varphi)^{\otimes 2}] -2 & \text{if $n=2$,}\\
1 & \text{if $n \in \{0,4\}$,}\\
0 & \text{otherwise.}\end{cases}
\]
\end{prop}
\begin{proof}
The cases $n=0,1$ and $n\geq 5$ are clear. The case $n=4$ follows from Lemma~\ref{lemm:det_tensor} and Lemma~\ref{lemm:det_symplectic}. The case $n=3$ then follows from the case $n=4$ and Lemma~\ref{lem:ext_n-1}. We now consider the case $n=2$. We have in $\GW^+(X)$
\begin{align*}
[\ext^2(E \otimes F,\varepsilon \otimes \varphi)]
&=\langle 2 \rangle [\Sym^2(E, \varepsilon)] + \langle 2 \rangle [\Sym^2(F,\varphi)]&\text{by \ref{lem:lambda_EF} and \ref{lemm:det_symplectic}} \\ 
&= [(E,\varepsilon)^{\otimes 2}] - \langle 2 \rangle + [(F,\varphi)^{\otimes 2}] - \langle 2 \rangle &\text{by \ref{lem:tens2_decomp} and \ref{lemm:det_symplectic}} 
\end{align*}
and $\langle 2 \rangle + \langle 2 \rangle =2 \in \GW^+(\Spec(\Zh))$, as evidenced by the computation
\[
\begin{pmatrix}
1&-1\\
1&1
\end{pmatrix}
\begin{pmatrix}
1&0\\
0&1
\end{pmatrix}
\begin{pmatrix}
1&1\\
-1&1
\end{pmatrix}
=
\begin{pmatrix}
2&0\\
0&2
\end{pmatrix}.
\qedhere
\]
\end{proof}

\section{Grothendieck--Witt groups of representations}
\label{sect:GW_rep}
Let $B$ be a commutative ring with $2\in B^\times$ and $G$ be a flat affine group scheme over $B$. Let $R_B$ be the abelian category of representations of $G$ over $B$, which are of finite type as $B$-modules. We let $P_B$ be the full subcategory of $R_B$ whose objects are projective as $B$-modules. The latter category is exact. If $P$ is an object of $P_B$, then its dual $P^\vee:=\Hom_B(P,B)$ is naturally endowed with an action of $G$ and thus can be seen as an object of $P_B$. The morphism of functors $\varpi\colon 1\simeq {}^{\vee\vee}$ is easily seen to be an isomorphism of functors $P_B\to P_B$, and it follows that $P_B$ is an exact category with duality. 

Let now $\Db(R_B)$, resp.\ $\Db(P_B)$, be the derived category of bounded complexes of objects of $R_B$, resp.\ $P_B$. The category $\Db(P_B)$ is a triangulated category with duality in the sense of Balmer (\cite[Definition 1.4.1]{Balmer05b}) and therefore one can consider its (derived) Witt groups $\W^i(\Db(P_B))$ (\cite[Definition 1.4.5]{Balmer05b}) that we denote by $\W^i(B; G)$ for simplicity. We can also consider the Grothendieck--Witt groups $\GW^i(\Db(P_B))$ (as defined in \cite[\S 2]{Walter03}) that we similarly denote by $\GW^i(B;G)$.

\begin{lem}\label{lem:Wittabelian}
Suppose that $B$ is a field of characteristic not two. For any $i\in\ZZ$, we have
\[
\W^{2i+1}(B;G)=0.
\]
\end{lem}
\begin{proof}
Since $P_B=R_B$, the category $\Db(R_B)$ is the derived category of an abelian category. We can thus apply \cite[Proposition 5.2]{Balmer02}.
\end{proof}

We now suppose that $A$ is a Dedekind domain with quotient field $K$ (we assume that $A \neq K$). We assume that $2 \in A^\times$, and let $G$ be a flat affine group scheme over $A$. Then we may consider the full subcategory $R_A^{\fl}$ of $R_A$ consisting of those representations of $G$ over $A$, which as $A$-modules are of finite length, or equivalently are torsion.

Any object of $\Db(P_A)$ has a well-defined support, and we can consider the (full) subcategory $\Dbfl(P_A)$ of $\Db(P_A)$ whose objects are supported on a finite number of closed points of $\Spec(A)$. This is a thick subcategory stable under the duality. As a consequence of \cite[Theorem~73]{Balmer05b}, we obtain a 12-term periodic long exact sequence
\begin{equation}\label{eq:12term}
\cdots\to \W^i(\Dbfl(P_A))\to \W^i(\Db(P_A))\to \W^i(\Db(P_A)/\Dbfl(P_A))\to \W^{i+1}(\Dbfl(P_A))\to \cdots
\end{equation}

We now identify the quotient category $\Db(P_A)/\Dbfl(P_A)$. Note that the extension of scalars induces a duality-preserving, triangulated functor $\Db(P_A)\to \Db(P_K)$ which is trivial on the subcategory $\Dbfl(P_A)$. (The category $\Db(P_K)$ is constructed by setting $B=K$ above, for the group scheme $G_K$ over $K$ obtained by base-change from $G$.) We thus obtain a duality-preserving, triangulated functor
\[
\Db(P_A)/\Dbfl(P_A)\to \Db(P_K).
\]

\begin{lem}
The functor $\Db(P_A)/\Dbfl(P_A)\to \Db(P_K)$ is an equivalence of triangulated categories with duality.
\end{lem}

\begin{proof}
We have a commutative diagram of functors
\[
\xymatrix{\Db(P_A)\ar[r]\ar[d]  & \Db(P_K)\ar[d] \\
\Db(R_A)\ar[r] & \Db(R_K)}
\]
in which the vertical arrows are equivalences (use \cite[\S 2.2, Corollaire]{Serre68b}. The composite $\Dbfl(P_A)\to \Db(P_A)\to \Db(R_A)$ has essential image the subcategory $\Dbfl(R_A)$ of objects of $\Db(R_A)$ whose homology is of finite length. As observed in \cite[Remarque, p.43]{Serre68b}, the functor $R_A \to R_K$ induces an equivalence $R_A/R_A^{\fl} \simeq R_K$.   Then it follows from \cite[\S 1.15, Lemma]{Keller99} that the induced functor $\Db(R_A)/\Dbfl(R_A) \to \Db(R_K)$ is an equivalence (the argument given in \cite[\S 1.15, Example (b)]{Keller99} works in the equivariant setting). The statement follows.
\end{proof}

As a consequence, the exact sequence \eqref{eq:12term} becomes
\begin{equation}
\label{eq:les_W}
\cdots\to \W^i(\Dbfl(P_A))\to \W^i(A;G)\to \W^i(K;G_K)\to \W^{i+1}(\Dbfl(P_A))\to \cdots
\end{equation}

Now, suppose that $M$ is a representation of $G$ over $A$ that is of finite length. By \cite[\S2.2, Corollaire]{Serre68b}, we have an exact sequence of representations
\begin{equation}
\label{eq:res}
0\to P_1\to P_0\to M\to 0
\end{equation}
where $P_0,P_1 \in P_A$. Note that the $A$-module $M$ is torsion, hence $M^\vee=\Hom_A(M,A)$ vanishes. We obtain an exact sequence, by dualizing
\[
0\to P_0^\vee\to P_1^\vee\to \Ext^1_{A}(M,A)\to 0
\]
and it follows that $M^{\sharp}:=\Ext^1_{A}(M,A)$ is naturally endowed with a structure of a representation of $G$ over $A$. The isomorphisms $P_0 \to (P_0^\vee)^\vee$ and $P_1 \to (P_1^\vee)^\vee$ induce an isomorphism $M \to (M^{\sharp})^{\sharp}$, which does not depend on the choice of the resolution \eqref{eq:res}. The association $M \mapsto M^{\sharp}$ in fact defines a duality on the category $R_A^{\fl}$. 

\begin{lem}
\label{lem:W_shift}
For every $i \in \ZZ$, there exists an isomorphism
\[
\W^{i+1}(\Dbfl(P_A)) \simeq \W^i(\Db(R^{\fl}_A)).
\]
\end{lem}
\begin{proof}
This follows from the existence of an equivalence of triangulated categories $\Dbfl(P_A) \to \Db(R^{\fl}_A)$, which is compatible with the duality $\sharp$ of $\Db(R^{\fl}_A)$, and the duality $\vee$ of $\Dbfl(P_A)$ shifted by $1$. This equivalence is constructed using word-for-word the proof of \cite[Lemma~6.4]{Balmer02}, where the categories $\VB_{\OO}, \OO
\text{-}\!\modul, \OO\text{-}\fl\text{-}\modul$ are replaced by $P_A,R_A,R_A^{\fl}$.
\end{proof}

\begin{lem}
\label{lem:W_even_vanish}
For every $i\in \ZZ$, we have $\W^{2i}(\Dbfl(P_A))=0$.
\end{lem}
\begin{proof}
In view of Lemma~\ref{lem:W_shift}, this follows from \cite[Proposition 5.2]{Balmer02}, as the category $R^{\fl}_A$ is abelian. 
\end{proof}

\begin{prop}
\label{prop:W_G_inj}
Let $A$ be a Dedekind domain with quotient field $K$, such that $2 \in A^\times$, and let $G$ be a flat affine group scheme over $A$. Then for every $i\in \ZZ$, the morphism $\W^{2i}(A;G)\to \W^{2i}(K;G_K)$ is injective.
\end{prop}
\begin{proof}
This follows from Lemma~\ref{lem:W_even_vanish} and the sequence \eqref{eq:les_W}.
\end{proof}

\begin{thm}
\label{prop:GW_G_inj}
Let $A$ be a Dedekind domain with quotient field $K$, such that $2 \in A^\times$, and let $G$ be a split reductive group scheme over $A$. Then for every $i \in \ZZ$, the morphism $\GW^{2i}(A;G) \to \GW^{2i}(K;G_K)$ is injective.
\end{thm}
\begin{proof}
We have a commutative diagram where rows are exact sequences (constructed in \cite[Theorem 2.6]{Walter03})
\[
\xymatrix{
\K_0(A;G)\ar[r]\ar[d] &  \GW^{2i-1}(A;G)\ar[r]\ar[d] &  \W^{2i-1}(A;G)\ar[r]\ar[d] & 0 \\
\K_0(K;G_K)\ar[r] &  \GW^{2i-1}(K;G_K)\ar[r] &  \W^{2i-1}(K;G_K)\ar[r] & 0}
\]
in which the vertical arrows are induced by the extension of scalars, and $K_0(A;G)$ (resp.\ $K_0(K;G_K)$) denotes the Grothendieck group of the triangulated category $\Db(P_A)$ (resp.\ $\Db(P_K)$). Denoting by $K_0(R_A)$ (resp.\ $K_0(R_K)$) the Grothendieck group of the category $R_A$ (resp.\ $R_K$), the natural morphisms $K_0(R_A) \to K_0(A;G)$ and $K_0(R_K) \to K_0(K;G)$ are isomorphisms (their inverses are constructed using the Euler characteristic). Since the morphism $\K_0(R_A) \to \K_0(R_K)$ is an isomorphism by \cite[Th\'eor\`eme 5]{Serre68b}, so is $K_0(A;G) \to K_0(K;G_K)$. On the other hand, we have $\W^{2i-1}(K;G_K)=0$ by Lemma~\ref{lem:Wittabelian}. We deduce that the morphism $\GW^{2i-1}(A;G) \to \GW^{2i-1}(K;G_K)$ is surjective.

Next consider the commutative diagram where rows are exact sequences (see again \cite[Theorem 2.6]{Walter03})
\[
\xymatrix{\GW^{2i-1}(A;G)\ar[r]\ar@{->>}[d]& \K_0(A;G)\ar[r]\ar[d]^{\simeq} &  \GW^{2i}(A;G)\ar[r]\ar[d] &  \W^{2i}(A;G)\ar[r]\ar@{^{(}->}[d] & 0 \\
\GW^{2i-1}(K;G_K)\ar[r] &\K_0(K;G_K)\ar[r] & \GW^{2i}(K;G_K)\ar[r] &  \W^{2i}(K;G_K)\ar[r] & 0}
\]
The indicated surjectivity and bijectivity have been obtained above, and the injectivity in Proposition~\ref{prop:W_G_inj}. The statement then follows from a diagram chase.
\end{proof}

\section{The \texorpdfstring{$\lambda$}{lambda}-operations}
\label{sect:lambda}
\numberwithin{thm}{subsection}
\numberwithin{equation}{subsection}
\renewcommand{\theequation}{\thesubsection.\alph{equation}}

Let $X$ be a scheme, and $G$ a flat affine group scheme over $X$. We denote by $\GW^+(X;G)$ and $\GW^-(X;G)$ the Grothendieck--Witt groups of the exact category of $G$-equivariant vector bundles over $X$. We set $\GW^{\pm}(X;G) = \GW^+(X;G) \oplus \GW^-(X;G)$. When $A$ is a commutative noetherian $\Zh$-algebra and $X =\Spec(A)$, by \cite[Theorem~6.1]{Walter03} we have natural isomorphisms $\GW^+(\Spec(A);G) \simeq \GW^0(A;G)$ and $\GW^-(\Spec(A);G) \simeq \GW^2(A;G)$ (in the notation of \S\ref{sect:GW_rep}).

\subsection{Exterior powers of metabolic forms}
Let $X$ be a scheme and $G$ a flat affine group scheme over $X$. Let $E \to X$ be a $G$-equivariant vector bundle. For $\varepsilon \in \{1,-1\}$, the associated hyperbolic $\varepsilon$-symmetric $G$-equivariant bundle over $X$ is
\[
H_\varepsilon(E) = \Big(E \oplus E^\vee,
\begin{pmatrix}
0 & 1\\
\varepsilon \varpi_E & 0
\end{pmatrix}\Big)
\]
where $\varpi_E \colon E \to (E^{\vee})^\vee$ is the canonical isomorphism. These constructions induce morphisms of abelian groups (see e.g.\ \cite[Proposition~2.2 (c), Theorem~6.1]{Walter03})
\begin{equation}
\label{eq:h_K_GW}
h_+ \colon K_0(X;G) \to \GW^+(X;G),\quad h_- \colon K_0(X;G) \to \GW^-(X;G)
\end{equation}
where $K_0(X;G)$ denotes the Grothendieck group of $G$-equivariant vector bundles on $X$.

\begin{lem}
\label{lem:ext_metabolic}
Let $M$ be a $G$-equivariant vector bundle over $X$ equipped with a $G$-equivariant nondegenerate $\varepsilon$-symmetric bilinear form $\mu$, for some $\varepsilon \in \{1,-1\}$. Assume that $(M,\mu)$ admits a ($G$-invariant) Lagrangian $L$, and let $n \in \NN$.
\begin{enumerate}[label={(\roman*)}]
\item
\label{lem:ext_metabolic:1}
The class $[\ext^n(M,\mu)] \in \GW^{\pm}(X;G)$ depends only on $n,\varepsilon$ and the $G$-equivariant vector bundle $L$ over $X$ (but not on $(M,\mu)$).

\item 
\label{lem:ext_metabolic:2}
If $n$ is odd, the $G$-equivariant nondegenerate $\varepsilon^n$-symmetric bilinear form $\ext^n(M,\mu)$ is metabolic.
\end{enumerate}
\end{lem}
\begin{proof}
We may assume that $X$ is connected. Let $Q=M/L$, and recall that $\mu$ induces an isomorphism $\varphi \colon Q \xrightarrow{\sim} L^\vee$. The vector bundle $\ext^nM$ is equipped with a decreasing filtration by $G$-invariant subsheaves
\[
(\ext^nM)^i = \im(\ext^iL \otimes \ext^{n-i}M \to \ext^nM)
\]
fitting into commutative squares
\begin{equation}
\label{squ:filtration}
\begin{gathered}
\xymatrix{
\ext^iL \otimes \ext^{n-i}M\ar[r] \ar[d] & (\ext^nM)^i  \ar[d] \\ 
\ext^iL \otimes \ext^{n-i}Q \ar[r] & (\ext^nM)^i/(\ext^nM)^{i+1}
}
\end{gathered}
\end{equation}
where the bottom horizontal arrow is an isomorphism (see e.g.\ \cite[V, Lemme 2.2.1]{SGA6}). Since $Q \simeq L^\vee$, this yields exact sequences of $G$-equivariant sheaves
\begin{equation}
\label{eq:ex_lambda_M_i}
0 \to (\ext^nM)^{i+1} \to (\ext^nM)^i \to \ext^iL \otimes \ext^{n-i}L^\vee \to 0
\end{equation}
from which we deduce by induction on $i$ that $(\ext^nM)^i$ is a subbundle of $\ext^nM$ (i.e.\ the quotient $\ext^nM/(\ext^nM)^i$ is a vector bundle). Assuming that $L$ has rank $r$, then $\ext^iL \otimes \ext^{n-i}L^\vee$ has rank $\binom{r}{i}\binom{r}{n-i}$. By induction on $i$, using the sequences \eqref{eq:ex_lambda_M_i}, we obtain
\[
\rank (\ext^nM)^i = \sum_{j=i}^r \binom{r}{j}\binom{r}{n-j}.
\]
An elementary computation with binomial coefficients then yields:
\begin{equation}
\label{eq:rank}
\rank (\ext^nM)^i + \rank (\ext^nM)^{n+1-i} = \rank \ext^nM.
\end{equation}
Let $i,j$ be integers. We have a commutative diagram
\[ \xymatrix{
\ext^iL \otimes \ext^{n-i}M \ar@{->>}[r] \ar[d]_{\alpha} &(\ext^nM)^i \ar@{^{(}->}[r] \ar[d] & \ext^nM \ar[d]^{\ext^n\mu}\\ 
(\ext^jL \otimes \ext^{n-j}M)^\vee& ((\ext^nM)^j)^\vee \ar@{_{(}->}[l]&  (\ext^nM)^\vee \ar@{->>}[l]
}\]
where $\alpha$ is defined by setting, for every open subscheme $U$ of $X$ and $x_1,\dots,x_i,y_1,\dots,y_j \in H^0(U,L)$ and $x_{i+1},\dots,x_n,y_{j+1},\dots,y_n \in H^0(U,M)$ (see \eqref{eq:explicit_ext:det})
\begin{equation}
\label{eq:alpha}
\alpha(x_1 \wedge \dots \wedge x_i \otimes x_{i+1} \wedge \dots \wedge x_n)(y_1 \wedge \dots \wedge y_j \otimes y_{j+1} \wedge \dots \wedge y_n) = \det(\mu(x_i,y_j)).
\end{equation}
If $i+j > n$, then for each $\sigma \in \mathfrak{S}_n$ there exists $e \in \{1,\dots,n\}$ such that $x_e \in H^0(U,L)$ and $y_{\sigma(e)} \in H^0(U,L)$, so that $\mu(x_e,y_{\sigma(e)})=0$, which by \eqref{eq:alpha}
 implies that $\alpha =0$. Thus $(\ext^nM)^i \subset ((\ext^nM)^j)^\perp$ in this case. In particular $(\ext^nM)^i$ is a sub-Lagrangian of $\ext^n(M,\mu)$ when $2i > n$.

If $n=2k-1$ with $k \in \NN$, then $2 \rank (\ext^nM)^k = \rank \ext^nM$ by \eqref{eq:rank}, hence the subbundle $(\ext^nM)^k$ is a Lagrangian in $\ext^n(M,\mu)$. This proves \ref{lem:ext_metabolic:2}. Moreover, it follows that the class of $\ext^n(M,\mu)$ in $\GW^\pm(X;G)$ coincides with the class of the hyperbolic form $H_{\varepsilon}((\ext^nM)^k)$, hence depends only on the class in $K_0(X;G)$ of the $G$-equivariant vector bundle $(\ext^nM)^k$ (see \eqref{eq:h_K_GW}). In view of the sequences \eqref{eq:ex_lambda_M_i}, the latter depends only on the classes of $\ext^iL \otimes \ext^{n-i}L^\vee$ in $K_0(X;G)$ for $i \geq k$, from which \ref{lem:ext_metabolic:1} follows when $n$ is odd.

Assume now that $n=2k$  with $k \in \NN$. Then the inclusion of the subbundle $(\ext^nM)^k \subset ((\ext^nM)^{k+1})^\perp$ is an equality by rank reasons (see \eqref{eq:rank}). By \cite[Proposition~2.2 (d), Theorem~6.1]{Walter03} we have
\begin{equation}
\label{eq:sublagrangian_red}
[(M,\mu)] = [H_1((\ext^nM)^{k+1})] + [((\ext^nM)^k/(\ext^nM)^{k+1},\rho)] \in \GW^+(X;G),
\end{equation}
where $\rho$ is the bilinear form induced by $\mu$ on $(\ext^nM)^k/(\ext^nM)^{k+1}$, which in view of \eqref{squ:filtration} is $G$-equivariantly isometric to the form $\beta$ fitting into the commutative diagram
\[ \xymatrix{
\ext^kL \otimes \ext^kM\ar@{->>}[r] \ar[d]^{\alpha} & \ext^k L \otimes \ext^kQ \ar[d] \ar[r]^{\simeq} & \ext^k L \otimes \ext^kL^\vee\ar[d]^\beta\\ 
(\ext^kL \otimes \ext^kM)^\vee & (\ext^k L \otimes \ext^kQ)^\vee\ar@{_{(}->}[l] & (\ext^k L \otimes \ext^kL^\vee)^\vee \ar[l]_{\simeq}
}\]
where the horizontal isomorphisms are induced by $\varphi \colon Q \xrightarrow{\sim} L^\vee$. The formula \eqref{eq:alpha} (and the fact that $\varphi$ is induced by $\mu$) yields the formula, for every open subscheme $U$ of $X$ and $x_1,\dots,x_k,y_1,\dots,y_k \in H^0(U,L)$ and $f_1,\dots,f_k,g_1,\dots,g_k \in H^0(U,L^\vee)$,
\[
\beta(x_1 \wedge \dots \wedge x_k \otimes f_1 \wedge \dots \wedge f_k,y_1 \wedge \dots \wedge y_k \otimes g_1 \wedge \dots \wedge g_k) = \det \begin{pmatrix}
0 & (\varepsilon g_j(x_i))\\
(f_i(y_j)) & 0
\end{pmatrix}
\]
(where $i,j$ run over $1,\dots,k$, and so the indicated determinant is $n \times n$), which shows that the bilinear form $\beta$ depends only on the $G$-equivariant vector bundle $L$ (and not on $\mu$). It follows that the isometry class of that $G$-equivariant form $((\ext^nM)^k/(\ext^nM)^{k+1},\rho)$ depends only on $L$. As above, the class of the hyperbolic form $H_1((\ext^nM)^{k+1})$ in $\GW^+(X;G)$ also depends only on $L$, so that \ref{lem:ext_metabolic:1} follows from \eqref{eq:sublagrangian_red} when $n$ is even.
\end{proof}

\subsection{The \texorpdfstring{$\lambda$}{lambda}-ring structure}
We will use the notion of (pre-)$\lambda$-rings, recalled in Appendix~\ref{appendix:lambda} below.

\begin{prop}
Let $X$ be a scheme and $G$ a flat affine group scheme over $X$. Then the exterior powers operations
\[
\lambda^i\colon \GW^{\pm}(X;G)\to \GW^{\pm}(X;G)
\]
defined by $(P,\varphi)\mapsto (\ext^iP,\ext^i\varphi)$ endow the ring $\GW^{\pm}(X;G)$ with the structure of a pre-$\lambda$-ring.
\end{prop}
\begin{proof}
The structure of the proof is the same as that of \cite[Proposition~2.1]{Zibrowius15}, and is based on the description of $\GW^{\pm}(X;G)$ in terms of generators and relations (see e.g.\ \cite[p.20]{Walter03}). It is clear that the exterior power operations descend to the set of isometry classes, and moreover the total exterior power operation is additive in the sense of \eqref{eq:ext_sum}. Finally, let $M$ is a $G$-equivariant vector bundle over $X$ equipped with a $G$-equivariant nondegenerate $\varepsilon$-symmetric bilinear form $\mu$, for some $\varepsilon \in \{1,-1\}$. If $(M,\mu)$ admits a $G$-equivariant Lagrangian $L$, then $L$ is also a $G$-equivariant Lagrangian in the hyperbolic form $H_{\varepsilon}(L)$, so that by Lemma~\ref{lem:ext_metabolic}~\ref{lem:ext_metabolic:1} the forms $\ext^n(M,\mu)$ and $\ext^n(H_{\varepsilon}(L))$ have the same class in $\GW^{\pm}(X;G)$.
%Step 1 and Step 2 are clear, while Step 3 follows from Lemma~\ref{lem:ext_metabolic}: indeed in the notation of that lemma, $L$ is also a Lagrangian in the hyperbolic form $H_{\varepsilon}(L)$, so that by \ref{lem:ext_metabolic:1} the forms $\ext^n(M,\mu)$ and $\ext^n(H_{\varepsilon}(L))$ have the same class in $\GW^{\pm}(X;G)$.
\end{proof}

\begin{prop}
\label{prop:GW0_lambda_ring}
Let $G$ be a split reductive group scheme over $\Zh$. Then the pre-$\lambda$-ring $\GW^+(\Spec(\Zh);G)$ is a $\lambda$-ring.
\end{prop}
\begin{proof}
By \cite[Proposition~2.1]{Zibrowius15} the pre-$\lambda$-ring $\GW^+(\Spec(\QQ);G_\QQ)$ is a $\lambda$-ring. It follows from Theorem~\ref{prop:GW_G_inj} that $\GW^+(\Spec(\Zh);G)$ is a pre-$\lambda$-subring of $\GW^+(\Spec(\QQ);G_\QQ)$, hence a $\lambda$-ring.
\end{proof}

\begin{cor}
\label{cor:GW0_lambda_ring}
For every $\Zh$-scheme $X$, the pre-$\lambda$-ring $\GW^+(X)$ is a $\lambda$-ring.
\end{cor}
\begin{proof}
This follows from Proposition \ref{prop:GW0_lambda_ring} (applied to the split reductive groups $\mathrm{O}_n$ and $\mathrm{O}_m \times \mathrm{O}_n$), using the arguments of \cite[Expos\'e VI, Th\'eor\`eme 3.3]{SGA6} (see \cite[\S3.2]{Zibrowius15} for details).
\end{proof}

When $x\in \GW^-(X)$ is the class of a rank two symplectic bundle, it follows from Lemma~\ref{lemm:det_symplectic} that $\lambda_t(x) = 1 +tx +t^2$ (see \eqref{def:lambda_t}). In other words, in the notation of \eqref{eq:ell}, we have
\begin{equation}
\label{eq:lambda_rank_2}
\lambda^i(x) = \ell_i(x) \in \GW^{\pm}(X) \quad \text{for all $i \in \NN\smallsetminus\{0\}$}.
\end{equation}

\begin{lem}
\label{lem:GW2_lambda_ring}
The relations \eqref{eq:L1} and \eqref{eq:L2} are satisfied for all $x,y,z \in \GW^-(X)$.
\end{lem}
\begin{proof}
By the symplectic splitting principle \cite[\S10]{Panin10pred}, we may assume that $x,y,z$ are each represented by a rank two symplectic bundle. In view of \eqref{eq:lambda_rank_2}, the relation \eqref{eq:L2} follows from Lemma \ref{lem:RZ}. The relation \eqref{eq:L1} has been verified in Proposition \ref{prop:lambda_22}, see Lemma \ref{lem:RXY}.
\end{proof}

\begin{thm}
\label{prop:GW_pm_lambda_ring}
For every $\Zh$-scheme $X$, the pre-$\lambda$-ring $\GW^{\pm}(X)$ is a $\lambda$-ring.
\end{thm}
\begin{proof}
Taking Proposition \ref{prop:GW0_lambda_ring} and Lemma \ref{lem:GW2_lambda_ring} into account, it only remains to verify \eqref{eq:L1} when $x \in \GW^+(X)$ and $y\in \GW^-(X)$. Let $i \geq n$, and consider the scheme $X \times \HP^i$. It is endowed with a universal symplectic bundle of rank two, whose class we denote by $u \in \GW^-(X \times \HP^i)$. Denote again by $x,y \in \GW^\pm(X \times \HP^i)$ the pullbacks of $x,y \in \GW^\pm(X)$. Then using successively Proposition \ref{prop:GW0_lambda_ring} and Lemma \ref{lem:GW2_lambda_ring}
\[
\lambda_t(xyu) = \lambda_t(x) \lambda_t(yu) = \lambda_t(x) \lambda_t(y) \lambda_t(u).
\]
On the other hand, by Lemma \ref{lem:GW2_lambda_ring}
\[
\lambda_t(xyu) = \lambda_t(xy) \lambda_t(u).
\]
The quaternionic projective bundle theorem \cite[Theorem~8.1]{Panin10pred} implies that the $\GWe(X)$-module $\GWe(X \times \HP^i)$ is free on the basis $1,u,\dots,u^i$. Modding out $\gamma-1$, we obtain a decomposition
\[
\GW^\pm(X \times \HP^i) = \GW^\pm(X) \oplus \GW^\pm(X)u \oplus \cdots \oplus \GW^\pm(X)u^i.
\]
In view of \eqref{eq:lambda_rank_2}, it follows from Lemma \ref{lem:RB} that the $u^n$-component of the $t^n$-coefficient of $\lambda_t(xy) \lambda_t(u)$ is $\lambda^n(xy)$, and that the $u^n$-component of the $t^n$-coefficient of $\lambda_t(x) \lambda_t(y) \lambda_t(u)$ is $P_n(\lambda^1(x),\dots,\lambda^n(x),\lambda^1(y),\dots,\lambda^n(y))$. This proves \eqref{eq:L1}.
\end{proof}

Let $X$ be a $\Zh$-scheme. In view of Lemma \ref{lem:lambdapowerseries2}, the $\lambda$-ring structure on $\GW^{\pm}(X)$ induces a $\lambda$-ring structure on $\widehat{\GW^{\pm}(X)} \simeq \GWe(X)$. Explicitly, denoting by $\rho \colon \GW^{\pm}(X) \to \GWe(X)$ the canonical homomorphism of abelian groups (see Appendix~\ref{sec:gradedrings}), we have for $i\in \ZZ$ and $n\in \NN$,
\begin{equation}\label{eqn:periodicitylambda}
\lambda^n(\rho(r)\cdot \gamma^i) = \begin{cases} \rho(\lambda^n(r)) \cdot  \gamma^{ni} & \text{ if $r\in \GW^+(X)$,} \\
\rho(\lambda^n(r)) \cdot\gamma^{\frac {n(2i+1)}2} & \text{ if $r\in \GW^-(X)$ and $n$ is even,} \\
\rho(\lambda^n(r)) \cdot  \gamma^{\frac {n(2i+1)-1}2}& \text{ if $r\in \GW^-(X)$ and $n$ is odd.} \end{cases}
\end{equation}

\section{The Adams operations}\label{sec:stable}
\subsection{The unstable Adams operations}
\label{sec:unstable}
The $\lambda$-operations constructed in \S\ref{sect:lambda} are not additive (with the exception of $\lambda^1$), and there is a standard procedure to obtain additive operations from the $\lambda$-operations which is valid in any pre-$\lambda$-ring, see e.g.\ \cite[\S5]{Atiyah-Tall}. Indeed, for any $\Zh$-scheme $X$, we define the (unstable) Adams operations
\[
\psi^n \colon \GW^{2i}_0(X) \to \GW^{2ni}_0(X) \quad \text{ for $n \in \NN\smallsetminus\{0\}$, $i \in \ZZ$}
\]
through the inductive formula (see e.g.\ \cite[Proof of Proposition~5.4]{Atiyah-Tall})
\begin{equation}
\label{eq:Adams_inductive}
\psi^{n}-\lambda^1\psi^{n-1}+\lambda^2{\psi}^{n-2}+\cdots+(-1)^{n-1}\lambda^{n-1}{\psi}^1+(-1)^nn\lambda^n=0.
\end{equation}
For instance, this yields
\[
\psi^1=\id \quad \text{ and } \quad \psi^2=(\id)^2-2\lambda^2.
\]
We also define $\psi^0$ as the composite
\begin{equation}
\label{eq:psi_0}
\psi^0 \colon \GW^{2i}_0(X) \xrightarrow{\rank} \ZZ^{\pi_0(X)} \to \GW^0_0(X).
\end{equation}

Assume now that $(E,\nu)$ is a rank two symplectic bundle on $X$, and let $x =[(E,\nu)] \in \GW^2_0(X)$ be its class. Then, by Lemma~\ref{lemm:det_symplectic}, we have for $n\in \NN \smallsetminus \{0\}$
\[
\lambda^n(x)=
\begin{cases}
x & \text{if $n=1$},\\
\gamma & \text{if $n=2$},\\
0 & \text{if $n \not \in \{1,2\}$}.
\end{cases}
\]
Thus \eqref{eq:Adams_inductive} yields the inductive formula for $x$ as above (the class of a rank two symplectic bundle)
\begin{equation}
\label{eq:rank_two}
\psi^n(x) = x \psi^{n-1}(x) - \gamma \psi^{n-2}(x) \quad \text{ for $n \geq 2$}.
\end{equation}

\begin{prop}
\label{prop:Adams_in_lambda_ring}
The operations $\psi^n \colon \GWe(X) \to \GWe(X)$ are ring morphisms for $n\in \NN$, and satisfy the relation $\psi^m \circ \psi^n = \psi^{mn}$ for $m,n \in \NN$.
\end{prop}
\begin{proof}
This follows from Theorem~\ref{prop:GW_pm_lambda_ring} (see for instance \cite[Propositions 5.1 and 5.2]{Atiyah-Tall}).
\end{proof}

\begin{rem}
The operations $\psi^n$ for $n <0$ are classically defined using duality; since by definition a nondegenerate symmetric (resp.\ skew-symmetric) form is isomorphic to its dual (resp.\ the opposite of its dual), in our situation we could set, for $n <0$
\[
\psi^n(x) =
\begin{cases}
\psi^{-n}(x) & \text{when $x \in \GW^{4i}_0(X)$ for $i\in \ZZ$,}\\
-\psi^{-n}(x) & \text{when $x\in \GW^{4i+2}_0(X)$ for $i \in \ZZ$,}
\end{cases}
\]
making Proposition~\ref{prop:Adams_in_lambda_ring} valid for $m,n \in \ZZ$.
\end{rem}

\subsection{Adams Operations on hyperbolic forms}

Let $X$ be a $\Zh$-scheme, and consider its Grothendieck group of vector bundles $K_0(X)$. The exterior power operations yield a $\lambda$-ring structure on $K_0(X)$ (and in particular Adams operations $\psi^n$ for $n \in \NN \smallsetminus \{0\}$, using the formula \eqref{eq:Adams_inductive}), such that the forgetful morphism
\begin{equation}
\label{eq:forgetful}
f\colon \GWe(X) \to K_0(X)
\end{equation}
(mapping $\gamma$ to $1$) is a morphism of $\lambda$-rings. In this section, we consider the hyperbolic morphisms $h_{2i} \colon K_0(X) \to \GWe(X)$ (defined just below). Those are of course not morphisms of $\lambda$-rings (not even ring morphisms), but as we will see in Proposition~\ref{prop:h_psi_odd}, they do satisfy some form of compatibility with the Adams operations.\\

We define morphisms
\begin{equation}
\label{eq:h_2i}
h_{2i} \colon K_0(X) \to \GW^{2i}_0(X) \quad \text{ for $i\in \ZZ$}
\end{equation}
by the requirements that $h_0=h_+$ and $h_2=h_-$ (see \eqref{eq:h_K_GW}) under the identifications $\GW_0^0(X) \simeq \GW^+(X)$ and $\GW^2_0(X) \simeq \GW^-(X)$, and for any vector bundle $E\to X$
\begin{equation}
\label{eq:h_periodic}
\gamma \cdot h_{2i}(E) = h_{2(i+2)}(E) \quad \text{ for $i\in \ZZ$}.
\end{equation}

\begin{lem}
\label{lem:proj_h}
Let $a \in K_0(X)$ and $b \in \GW^{2j}_0(X)$. Then, in the notation of \eqref{eq:forgetful} and \eqref{eq:h_2i}, we have for any $i \in \ZZ$
\[
h_{2i}(a) \cdot b = h_{2(i+j)}(a \cdot f(b)).
\]
\end{lem}
\begin{proof}
Let $\varepsilon,\varepsilon' \in \{1,-1\}$. Let us consider vector bundles $A,B$ on $X$, and a nondegenerate $\varepsilon$-symmetric bilinear form $\nu$ on $B$. The isomorphism
\[
(A \otimes B) \oplus (A^\vee \otimes B) \xrightarrow{\begin{pmatrix}1 & 0\\0 &1\otimes \nu  \end{pmatrix}} (A \otimes B) \oplus (A^\vee \otimes B^\vee) \simeq (A \otimes B) \oplus (A\otimes B)^\vee
\]
induces an isometry
\[
\Big((A \otimes B) \oplus (A^{\vee} \otimes B),
\begin{pmatrix}
0 & 1 \otimes \nu\\
\varepsilon' \varpi_A \otimes \nu &0
\end{pmatrix}
\Big)
\simeq
\Big((A \otimes B) \oplus (A \otimes B)^\vee,
\begin{pmatrix}
0 & 1\\
\varepsilon \varepsilon' \varpi_{A \otimes B} &0
\end{pmatrix}
\Big),
\]
as evidenced by the computation
\begin{align*}
\begin{pmatrix}
1&0\\
0 & 1 \otimes \nu^\vee
\end{pmatrix}
\begin{pmatrix}
0&1\\
\varepsilon \varepsilon' \varpi_A \otimes \varpi_B & 0
\end{pmatrix}
\begin{pmatrix}
1&0\\
0 & 1 \otimes \nu
\end{pmatrix}
& =  
\begin{pmatrix}
0&1 \otimes \nu\\
\varepsilon \varepsilon' \varpi_A \otimes (\nu^\vee \circ \varpi_B) & 0
\end{pmatrix}
 \\
& =  
\begin{pmatrix}
0&1 \otimes \nu\\
\varepsilon' \varpi_A \otimes \nu & 0
\end{pmatrix}.
\end{align*}
The lemma follows.
\end{proof}

\begin{lem}
\label{lem:product_h}
For any $i,j\in \ZZ$ we have $h_{2i}(1)h_{2j}(1) = 2 h_{2(i+j)}(1)$.
\end{lem}
\begin{proof}
Take $a=1 \in K_0(X)$ and $b = h_{2j}(1) \in \GW^{2j}_0(X)$ in Lemma~\ref{lem:proj_h}.
\end{proof}

Observe that the classes $h$ and $\tau$ (see Notation~\ref{not:tau}) coincide respectively with $h_0(1)$ and $h_2(1)$. Thus Lemma~\ref{lem:product_h} implies that
\begin{equation}
\label{eq:tau_sq}
h^2=2h \quad ; \quad h\tau = 2\tau \quad ; \quad \tau^2 = 2\gamma h.
\end{equation}
Combining the relations $h\tau = 2\tau$ and $h =1 + \langle -1 \rangle$ yields
\begin{equation}
\label{eq:2_sigma:2}
\langle -1 \rangle\tau = \tau.
\end{equation}

\begin{lem}
\label{lem:psi_tau}
For $n \in \NN$, we have in $\GW^{2n}_0(\Spec(\Zh))$ (see Notation~\ref{not:tau})
\[
\psi^n(\tau)=\begin{cases} \tau\gamma^{\frac {n-1}2} & \text{ if $n$ is odd.} \\ 2 \langle -1 \rangle^{\frac{n}{2}}\gamma^{\frac{n}{2}} & \text{ if $n$ is even.} 
\end{cases}
\]
\end{lem}
\begin{proof}
We prove the lemma by induction on $n$, the cases $n=0,1$ being clear. If $n \geq 2$, we have by \eqref{eq:rank_two}
\begin{equation}
\label{eq:ind_hyperbolic:2}
\psi^n(\tau) = \tau \psi^{n-1}(\tau) -\gamma \psi^{n-2}(\tau).
\end{equation}
Assume that $n$ is odd. Using the induction hypothesis together with \eqref{eq:2_sigma:2} we obtain
\[
\tau \psi^{n-1}(\tau) = 2 \langle -1 \rangle^{\frac{n-1}{2}}\tau \gamma^{\frac{n-1}{2}} = 2 \tau \gamma^{\frac{n-1}{2}}.
\]
On the other hand, by induction we have
\[
\gamma \psi^{n-2}(\tau) = \gamma \tau \gamma^{\frac{n-3}{2}} = \tau \gamma^{\frac{n-1}{2}}.
\]
Combining these two computations with \eqref{eq:ind_hyperbolic:2} proves the statement when $n$ is odd. 

Assume now that $n$ is even. Using the induction hypothesis we have
\[
\tau \psi^{n-1}(\tau) = \tau^2 \gamma^{\frac{n-2}{2}} = 2h \gamma^{\frac{n}{2}} = 2(\langle -1 \rangle^{\frac{n}{2}} + \langle -1 \rangle^{\frac{n-2}{2}}) \gamma^{\frac{n}{2}},
\]
as well as
\[
\gamma \psi^{n-2}(\tau) = 2 \langle -1 \rangle^{\frac{n-2}{2}}\gamma^{\frac{n}{2}},
\]
and the result follows as above from \eqref{eq:ind_hyperbolic:2} when $n$ is even.
\end{proof}

\begin{prop}
\label{prop:h_psi_odd}
Let $E \to X$ be a vector bundle, and $n \in \NN, i \in \ZZ$. For $j \in \ZZ$, let us denote by $I_j$ the image of $h_{2j} \colon K_0(X) \to \GW^{2j}_0(X)$. 
\begin{enumerate}[label={(\roman*)}]
\item \label{prop:h_psi_odd:1}
If $n$ is odd, then $\lambda^n\circ h_{2i}(E)$ lies in $I_{in}$.

\item \label{prop:h_psi_odd:2}
If $n$ is odd, then $\psi^n \circ h_{2i}(E)$ lies in $I_{in}$.
\item \label{prop:h_psi_odd:3}
If $n$ is even, then $\psi^n \circ h_{2i}(E)$ lies in $2\GW_0^{2in}(X) + I_{in}$.
\end{enumerate}
\end{prop}
\begin{proof}
Statement \ref{prop:h_psi_odd:1} follows from Lemma~\ref{lem:ext_metabolic} \ref{lem:ext_metabolic:2} with $G=1$ (observe that by construction of the Grothendieck--Witt group, the classes of metabolic forms belong to the subgroup $I_{in} \subset \GW^{2in}_0(X)$). Let us prove \ref{prop:h_psi_odd:2} by induction on $n$. This is clear when $n=1$. Assume that $n$ is odd. When $j\in \{1,\dots,n-1\}$ is even, the element $\psi^{n-j} \circ h_{2i}(E)$ belongs to $I_{i(n-j)}$ by induction. When $j\in \{1,\ldots,n\}$ is odd the element $\lambda^j \circ h_{2i}(E)$ belongs to $I_{ij}$ by \ref{prop:h_psi_odd:1}. Since $I_{ik}\cdot \GW^{2i(n-k)}_0(X) \subset I_{in}$ for all $k \in \ZZ$ by Lemma~\ref{lem:proj_h}, it follows from the inductive formula \eqref{eq:Adams_inductive} that $\psi^n \circ h_{2i}(E)$ belongs to $I_{in}$. The proof of \ref{prop:h_psi_odd:3} is similar, noting that $n\lambda^n \circ h_{2i}(E)$ is divisible by $2$ (the starting case $n=0$ being clear from \eqref{eq:psi_0}).
\end{proof}

Recall the exact sequence of \cite[Theorem 2.6]{Walter03}, for $i\in \ZZ$,
\[
K_0(X) \xrightarrow{h_{2i}} \GW^{2i}_0(X) \to \W^{2i}(X) \to 0.
\]
When $X \neq \varnothing$, the $\lambda$-ring structure on $\GWe(X)$ does not descend to its quotient $\bigoplus_{i\in\ZZ} \W^{2i}(X)$, for instance because $\lambda^2(h_0(1)) = \langle -1\rangle$ has nonzero image in the Witt ring. However, Proposition~\ref{prop:h_psi_odd} implies the following:
\begin{cor}
Let $n \in \NN$ be odd. Then the operations $\psi^n,\lambda^n \colon \GW^{2i}_0(X) \to \GW^{2in}_0(X)$ descend to operations
\[
\psi^n,\lambda^n \colon \W^{2i}(X) \to \W^{2ni}(X).
\]
\end{cor}

\begin{rem}
If $-1$ is a square in $H^0(X,\OO_X)$, then $2=h_0(1) \in \GW^0_0(X)$. Therefore Proposition~\ref{prop:h_psi_odd}~\ref{prop:h_psi_odd:3} implies that the operation $\psi^n$ does descend to the Witt groups when $n$ is even (even though $\lambda^n$ does not).
\end{rem}

\subsection{Adams operations on the universal rank two bundle}
In this section, we consider the universal symplectic bundle $(U,\varphi)$ over $\HP^1$, and denote by $u$ its class in $\GW^2_0(\HP^1)$. 

\begin{prop-df}
\label{prop:omega_n}
Let $n \in \NN$. There exists a unique element
\[
\omega(n) \in \GW^{2n-2}_0(\Spec(\Zh))
\]
such that
\[
\psi^n(u-\tau) = \omega(n)\cdot (u -\tau) \in \GW^{2n}_0(\HP^1).
\]
\end{prop-df}
\begin{proof}
The first Borel class of $(U,\varphi)$ in $\GW^2_0(\HP^1)$ is $u-\tau$ (see \cite[Theorem~9.9]{Panin10prec}). By the quaternionic projective bundle theorem \cite[Theorem~8.1]{Panin10pred}, the $\GWe(\Spec(\Zh))$-module $\GWe(\HP^1)$ is free on the basis $1,u-\tau$. This implies in particular the uniqueness part of the statement. Let us write
\[
\psi^n(u-\tau) = a + b (u-\tau)
\]
with $a\in  \GW^{2n}_0(\Spec(\Zh))$ and $b \in \GW^{2n-2}_0(\Spec(\Zh))$.
Consider the morphism of $\Zh$-schemes $i_0 \colon \Spec(\Zh)=\HP^0 \to \HP^1$ of \eqref{eq:i_n}. Since $i_0^*(u) = \tau$, we have
\[
a = i_0^*(a + b (u-\tau)) = i_0^*\circ \psi^n(u-\tau) = \psi^n \circ i_0^*(u-\tau)=\psi^n(0)=0.
\]
So we may set $\omega(n)=b$.
\end{proof}

\begin{lem}
\label{lem:omega_psi}
Let $m,n \in \NN$. Then $\omega(mn) = \omega(n) \cdot \psi^n(\omega(m))$.
\end{lem}
\begin{proof}
Indeed by Proposition~\ref{prop:Adams_in_lambda_ring}, we have in $\GW^{2mn}_0(\HP^1)$
\begin{align*}
\psi^{mn}(u-\tau)
&= \psi^n \circ \psi^m(u-\tau) \\ 
&= \psi^n(\omega(m) \cdot (u-\tau)) \\
&= \psi^n(\omega(m)) \cdot \psi^n(u-\tau) \\
&= \omega(n)\cdot \psi^n(\omega(m))\cdot (u-\tau). \qedhere
\end{align*}
\end{proof}

From the inductive definition of the Adams operations, we deduce an inductive formula for the classes $\omega(n)$:
\begin{lem}
\label{lem:inductive_omega_n}
We have $\omega(0)=0, \omega(1)=1$, and if $n\geq 2$
\[
\omega(n) = \tau \omega(n-1) -\gamma \omega(n-2) + \psi^{n-1}(\tau).
\]
\end{lem}
\begin{proof}
The computations of $\omega(0)$ and $\omega(1)$ are clear. Assume that $n \geq 2$. Then by \eqref{eq:rank_two} we have in $\GW^{2n}_0(\HP^1)$
\begin{eqnarray*}
\psi^n(u-\tau) & = & \psi^n(u)-\psi^n(\tau) \\
& = & u\psi^{n-1}(u)-\gamma\psi^{n-2}(u)-\tau\psi^{n-1}(\tau)+\gamma\psi^{n-2}(\tau) \\
& = & u\psi^{n-1}(u-\tau)+(u-\tau)\psi^{n-1}(\tau)-\gamma\psi^{n-2}(u-\tau).
\end{eqnarray*}

By the quaternionic projective bundle theorem \cite[Theorem~8.1]{Panin10pred} we have $(u-\tau)^2=0$, hence $u(u-\tau)=\tau(u-\tau)$, so that
\[
\psi^n(u-\tau) =(u-\tau)\Big( \tau \omega(n-1) + \psi^{n-1}(\tau) -\gamma \omega(n-2)\Big),
\]
from which the result follows.
\end{proof}

We are now in position to find an explicit expression for the elements $\omega(n)$. For this, recall from Notation~\ref{not:tau} that $h=1-\epsilon$.
\begin{prop}
\label{prop:omega_explicit}
We have
\[
\omega(n) = 
\begin{cases}
\displaystyle{n \Big(\frac{n-1}{2} h + \langle -1\rangle^{\frac{n-1}{2}}\Big)\gamma^{\frac{n-1}{2}} }& \text{if $n$ is odd},\\
\displaystyle{\frac{n^2}{2}\tau\gamma^{\frac{n-2}{2}}} & \text{if $n$ is even}.
\end{cases}
\]
\end{prop}
\begin{proof}
We proceed by induction on $n$, the cases $n=0,1$ being clear. Let $n \geq 2$. Assume that $n$ is even. Recall that  $h\tau = 2\tau$ by \eqref{eq:tau_sq} and that $\tau \langle -1 \rangle = \tau$ by \eqref{eq:2_sigma:2}. Combining these observations with the explicit formula for $\omega(n-1)$ (known by induction) yields
\[
\tau \omega(n-1) = (n-1)^2\tau \gamma^{\frac{n-2}{2}},
\]
hence, using the inductive hypothesis and Lemma~\ref{lem:psi_tau}
\[
\tau \omega(n-1) - \gamma \omega(n-2) + \psi^{n-1}(\tau) = (n-1)^2\tau \gamma^{\frac{n-2}{2}} - \frac{(n-2)^2}{2}\tau \gamma^{\frac{n-2}{2}} + \tau \gamma^{\frac{n-2}{2}} = \frac{n^2}{2} \tau \gamma^{\frac{n-2}{2}},
\]
which coincides with $\omega(n)$ by Lemma~\ref{lem:inductive_omega_n}, as required. 

Assume that $n$ is odd. Observe that $h=\langle -1 \rangle ^{\frac{n-1}{2}} + \langle -1 \rangle^{\frac{n-3}{2}}$, so that we have by induction
\[
\omega(n-2) = (n-2)\Big(\frac{n-1}{2}h- \langle -1 \rangle^{\frac{n-1}{2}}\Big) \gamma^{\frac{n-3}{2}} .
\]
Therefore, using Lemma~\ref{lem:inductive_omega_n}, Lemma~\ref{lem:psi_tau} and \eqref{eq:tau_sq} (and the inductive hypothesis)
\begin{align*}
\omega(n) 
&= \tau \omega(n-1) - \gamma \omega(n-2) + \psi^{n-1}(\tau) \\ 
&= \frac{(n-1)^2}{2} \tau^2 \gamma^{\frac{n-3}{2}} - (n-2) \Big(\frac{n-1}{2}h- \langle -1 \rangle^{\frac{n-1}{2}}\Big) \gamma^{\frac{n-1}{2}} + 2 \langle -1 \rangle^{\frac{n-1}{2}}\gamma^{\frac{n-1}{2}}\\
&= \Big((n-1)^2h-(n-2)\frac{n-1}{2}h + (n-2)\langle -1 \rangle^{\frac{n-1}{2}} + 2\langle -1 \rangle^{\frac{n-1}{2}}\Big)\gamma^{\frac{n-1}{2}} \\
&= n  \Big(\frac{n-1}{2}h+ \langle -1 \rangle^{\frac{n-1}{2}}\Big)\gamma^{\frac{n-1}{2}}.\qedhere
\end{align*}
\end{proof}

\subsection{Inverting \texorpdfstring{$\omega(n)$}{\textomega (n)}}
In order to define the stable Adams operations, we will be led to invert the elements $\omega(n) \in \GW^{2n-2}_0(\Spec(\Zh))$. Let us first observe that it is equivalent to invert somewhat simpler elements.

\begin{df}
For $n\in \NN$, we define an element $n^{\star} \in \GW^0_0(\Spec(\Zh))$ by
\[
n^{\star} = 
\begin{cases}
n & \text{ if $n$ is odd,}\\
\frac{n}{2}h & \text{ if $n$ is even}.
\end{cases}
\]
(Recall from Notation~\ref{not:tau} that $h=1-\epsilon \in \GW^0_0(\Spec(\Zh))$ is the hyperbolic class.)
\end{df}

\begin{lem}
\label{lemm:invert_omega_nstar}
Let $R =\GWe(\Spec(\Zh))$. Then the $R$-algebras $R[\frac{1}{n^\star}]$ and $R[\frac{1}{\omega(n)}]$ are isomorphic.
\end{lem}
\begin{proof}
We use the explicit formulas of Proposition~\ref{prop:omega_explicit}. Assume that $n$ is odd. Since $n=n^\star$ divides $\omega(n)$, it is invertible in $R[\frac{1}{\omega(n)}]$. Conversely, writing $n=2m+1$ we have (recall that $\epsilon = - \langle -1 \rangle$, so that $\epsilon^2 =1$)
\begin{align*}
\omega(n)\cdot (m (1+\epsilon) + \epsilon^m)
&= \gamma^mn(m(1-\epsilon) + (-\epsilon)^m) \cdot (m (1+\epsilon) + \epsilon^m) \\ 
&= \gamma^mn(m(1-\epsilon)\epsilon^m +m(1+\epsilon)(-\epsilon)^m + (-1)^m)\\
&= \gamma^mn(m\epsilon^m(1-\epsilon + (-1)^m(1+\epsilon)) + (-1)^m)\\
&= \gamma^mn(2m+1)(-1)^m=\gamma^mn^2(-1)^m
\end{align*}
(where the penultimate equality is seen for instance by distinguishing cases according to the parity of $m$). It follows that $\omega(n)$ is invertible in $R[\frac{1}{n^\star}] = R[\frac{1}{n}]$.

Now assume that $n$ is even. Then, by \eqref{eq:tau_sq}
\begin{equation}
\label{eq:omega_sq}
\omega(n)^2 =\Big(\frac{n^2}{2}\tau\Big)^2\gamma^{n-2}  = \frac{n^4}{2}h \gamma^{n-1} =n^3 n^\star \gamma^{n-1},
\end{equation}
so that $n^\star$ is invertible in $R[\frac{1}{\omega(n)}]$. On the other hand, using \eqref{eq:tau_sq}, we have
\[
(n^\star)^2 = \frac{n^2}{2}h = n \frac{n}{2}h,
\]
hence $n$ is invertible in $R[\frac{1}{n^\star}]$. Thus \eqref{eq:omega_sq} implies that $\omega(n)$ is invertible in $R[\frac{1}{n^\star}]$.
\end{proof}

We want now to formally invert the action of $n^{\star}$ on the spectrum $\sOGr$.
\begin{df}
We consider the ring
\begin{equation}
\label{def:B}
B=\ZZ[e]/(e^2-1),
\end{equation}
and for $n\in \NN$, we define an element $n^* \in B$ by
\[
n^* = 
\begin{cases}
n & \text{ if $n$ is odd,}\\
\frac{n}{2}(1-e) & \text{ if $n$ is even}.
\end{cases}
\]
\end{df}

For any $m,n \in \NN$, we have
\begin{equation}
\label{eq:m_n_*}
(mn)^* = m^*n^* \in B.
\end{equation}

\begin{rem}
Observe that the ring morphism $B \to \GWe(\Spec(\Zh))$ given by $e\mapsto \epsilon$ maps $n^*$ to $n^\star$. The ring $B$ may be identified with $\GW^+(\Spec \ZZ)$, but we will not use this observation.
\end{rem}

Denote by $\mathbb{S} \in \SH(\Zh)$ the sphere spectrum. Recall that each invertible element $u \in (\Zh)^\times$ defines an endomorphism $\langle u \rangle \in \mathrm{End}_{\SH(\Zh)}(\mathbb{S})$ (see e.g.\ \cite[2.2.8]{Deglise20}). Thus we may define a ring homomorphism
\begin{equation}
\label{eq:B_End_1}
B \to \mathrm{End}_{\SH(\Zh)}(\mathbb{S}), \quad e \mapsto -\langle -1 \rangle,
\end{equation}
which allows us to see $n^*$ as an endomorphism of the sphere spectrum and perform the formal inversion of $n^*$ in an efficient way as explained in \cite[\S 6]{Bachmann18b}. In short, we consider the diagram
\[
\mathbb{S}\xrightarrow{n^*}\mathbb{S}\xrightarrow{n^*}\ldots
\]
and define $\mathbb{S}[\frac {1}{n^*}]$ to be its homotopy colimit in $\SH(\Zh)$. Further, we set 
\[
\sOGr\Big[\frac{1}{n^*}\Big]:=\sOGr\wedge \phantom{i}\mathbb{S}\Big[\frac {1}{n^*}\Big].
\]
This is naturally a motivic ring spectrum. 

The $B$-algebra structure on $\GWe(\Spec (\Zh))$ induced by \eqref{eq:B_End_1} is given by $e \mapsto \epsilon$ (the argument is detailed in the last paragraph of the proof of \cite[Theorem~11.1.5]{Panin10prec}), and in particular maps $n^*$ to $n^{\star}$. It thus follows from Lemma~\ref{lemm:invert_omega_nstar}, that for any $i \in \NN$, the morphism $\sOGr[\frac{1}{n^*}] \to \Sigma_\Tc^{i(n-1)}\sOGr[\frac{1}{n^*}]$ induced by multiplication by $\omega(n)^i \in \GW^{2i(n-1)}_0(\Spec (\Zh))$ admits an inverse in $\SH(\Zh)$
\begin{equation}
\label{eq:omega_n_i}
\omega(n)^{-i} \colon \Sigma_\Tc^{i(n-1)}\sOGr\Big[\frac{1}{n^*}\Big] \to \sOGr\Big[\frac{1}{n^*}\Big].
\end{equation}

For any $i,j\in \ZZ$ and any smooth $\Zh$-scheme $X$, the spectrum $\Sigma_{S^1}^{j}\Sigma_{\PP^1}^i\Sigma_{\PP^1}^{\infty} X_+$ is a compact object in $\SH(\Zh)$ by \cite[\S 2.2, Lemma 2.2]{Jardine00} and it follows from \cite[Tag 094A]{stacks} that we have a canonical isomorphism
\begin{equation}
\label{eq:GW_inverted}
\Big[\Sigma_{\PP^1}^{\infty}X_+,\Sigma_{S^1}^{-j}\Sigma_{\PP^1}^i\Big(\sOGr\Big[\frac{1}{n^*}\Big]\Big)\Big]_{\SH(\Zh)}=\GW^{i}_j(X)\Big[\frac{1}{n^*}\Big]
\end{equation}
In case $X$ is merely a regular $\Zh$-scheme, the same property holds using the spectrum $p_X^*(\sOGr[\frac{1}{n^*}])$, where $p_X \colon X \to \Spec(\Zh)$ is the structural morphism.

\subsection{The stable Adams operations}
Recall from \cite[\S 8]{Panin10prec} and \cite[Theorem 1.3]{Schlichting13} that $\GW^2$ is naturally isomorphic to the object $\ZZ \times \HGr$ in the homotopy category $\Ho(\Zh)$, where $\HGr$ denotes the infinite quaternionic Grassmannian. Thus, by \cite[Theorem 4.1.4]{Deglise20} (where $\GW^2_0(X)$ is denoted $KSp_0(X)$) the Adams operations $\psi^n \colon \GW^2_0(X) \to \GW^{2n}_0(X)$ constructed in \S\ref{sec:unstable}, where $X$ runs over the smooth $\Zh$-schemes, are induced by a unique morphism $\psi^n:\GW^2\to \GW^{2n}$ in $\Ho(\Zh)$. Using the periodicity isomorphisms \eqref{eq:periodicity}, we obtain Adams operations:
\begin{equation}
\label{eq:psi_space}
\psi^n \colon \GW^{2i}\to \GW^{2ni}\; \text{ for $i$ odd}.
\end{equation}

We will need the following complement to \cite[Theorem 4.1.4]{Deglise20}:
\begin{lem}
\label{lemm:lift}
Let $E \in \SH(\Zh)$ be a $\Sp$-oriented ring spectrum, and consider for $a,b\in \ZZ$ the pointed motivic space $\Ec = \Omega^{\infty}_{\PP^1}\Sigma_{S^1}^a \Sigma_{\PP^1}^bE \in \Ho(\Zh)$. Let $i_1,\dots,i_r \in \ZZ$ be odd integers. Then each map $\GW^{2i_1} \wedge \dots \wedge \GW^{2i_r} \to \Ec$ in $\Ho(\Zh)$ is determined by the induced maps
\begin{equation}
\label{eq:lift_induced_maps}
\GW^{2i_1}_0(X_1) \times \dots \times \GW^{2i_r}_0(X_r) \to [(X_1 \times \dots \times X_r)_+,\Ec]_{\Ho(\Zh)}
\end{equation}
where $X_1,\dots,X_r$ run over the smooth $\Zh$-schemes.
\end{lem}
\begin{proof}
By \cite[\S 8]{Panin10prec}, for $j=1,\dots,r$, the pointed motivic space $\GW^{2i_j}$ can be expressed as a (homotopy) colimit of pointed smooth $\Zh$-schemes $Y_{m,j} = \{-m,\dots,m\} \times \mathrm{HGr}_{m,j}$ over $m \in \NN$, where $\mathrm{HGr}_{m,j}$ denotes an appropriate symplectic Grassmannian. Set $Y_m = Y_{m,1} \wedge \dots \wedge Y_{m,r}$ and $\Gc= \GW^{2i_1} \wedge \dots \wedge \GW^{2i_r}$. Then by \cite[Theorem~10.1]{Panin10prec} we have an exact sequence
\[
0 \to {\lim_m}^1 [S^1 \wedge Y_m,\Ec]_{\Ho(\Zh)} \to [\Gc,\Ec]_{\Ho(\Zh)} \to \lim_m [Y_m,\Ec]_{\Ho(\Zh)} \to 0.
\]
The $\lim^1$-term vanishes by \cite[Theorems~9.4,13.2,13.3]{Panin10prec} (see the proof of \cite[Theorem~13.1]{Panin10prec}). Thus a map $\Gc \to \Ec$ in $\Ho(\Zh)$ is determined by its restrictions to $[Y_m,\Ec]_{\Ho(\Zh)}$, for $m\in \NN$, each of which is determined by its restriction to $[(Y_{m,1} \times \dots \times Y_{m,r})_+,\Ec]_{\Ho(\Zh)}$, in view of \cite[Lemma~7.6]{Panin10prec}. The latter is the image of the tuple of canonical maps
\[
(Y_{m,1} \to \GW^{2i_1}, \dots , Y_{m,r} \to \GW^{2i_r}) \in \GW^{2i_1}_0(Y_{m,1}) \times \dots \times \GW^{2i_r}_0(Y_{m,r})
\]
under the map \eqref{eq:lift_induced_maps}.
\end{proof}

We are now in position to follow the procedure described in \cite[\S 4]{Deglise20} to construct the $n$-th stable Adams operation, for $n \in \NN$. For any integer $i\in \ZZ$, consider the motivic space
\[
\GWibig{2i}{n} = \Omega_{\Tc}^{\infty}\Sigma_\Tc^{i}\sOGr\Big[\frac{1}{n^*}\Big] \in \Ho(\Zh).
\]
The composite $\Sigma^\infty_{\Tc} \GW^{2i} \to \Sigma^i_{\Tc} \sOGr \to \Sigma^i_{\Tc} \sOGr[\frac{1}{n^*}]$ in $\SH(\Zh)$ yields by adjunction a morphism
\begin{equation}
\label{eq:GW_GWi}
\GW^{2i} \to \GWibig{2i}{n} \; \text{ in $\Ho(\Zh)$},
\end{equation}
while the morphism in $\SH(\Zh)$
\[
\Sigma_\Tc^\infty\Big(\Tc \wedge \GWibig{2i}{n}\Big) = \Sigma_\Tc \Sigma_\Tc^\infty \GWibig{2i}{n} \xrightarrow{\Sigma_\Tc(\text{counit})}\Sigma_\Tc\Sigma_\Tc^i \sOGr\Big[\frac{1}{n^*}\Big] = \Sigma_\Tc^{i+1} \sOGr\Big[\frac{1}{n^*}\Big]
\]
yields a morphism
\begin{equation}
\label{eq:GWi_bond}
\Tc \wedge \GWibig{2i}{n} \to \GWibig{2(i+1)}{n} \; \text{ in $\Ho(\Zh)$}.
\end{equation}
Using the morphism $\omega(n)^{-i}$ of \eqref{eq:omega_n_i}, we define a morphism in $\Ho(\Zh)$, for $i$ odd,
\begin{equation}
\label{eq:def_Psi}
\Psi^n_i \colon \GW^{2i} \xrightarrow{\psi^n} \GW^{2ni} \xrightarrow{\eqref{eq:GW_GWi}} \GWibig{2ni}{n} \xrightarrow{\Omega_{\Tc}^{\infty}\Sigma_{\Tc}^i\omega(n)^{-i}} \GWibig{2i}{n}.
\end{equation}

\begin{prop}
\label{prop:omegai}
Let $i\in\ZZ$ be odd and let $n \in \NN$. Then the diagram 
\[ \xymatrix{
\Tc^{\wedge 2} \wedge \GW^{2i}\ar[rr] \ar[d]_{\id_{\Tc^{\wedge 2}} \wedge \omega(n)^2 \psi^n} && \GW^{2(i+2)} \ar[d]^{\psi^n} \\ 
\Tc^{\wedge 2} \wedge \GW^{2n(i+2)-4} \ar[rr] && \GW^{2n(i+2)}
}\]
commutes in $\Ho(\Zh)$, where the horizontal arrows are induced by the bonding map $\sigma$ of \eqref{eq:def_sigma}.
\end{prop}
\begin{proof}
Let $X$ be a smooth $\Zh$-scheme, and denote by $p \colon \HP^1 \times \HP^1 \times X \to X$ the projection. Let $u_1,u_2 \in \GW^2_0(\HP^1 \times \HP^1 \times X)$ be the pullbacks of $u \in \GW^2_0(\HP^1)$ under the two projections. Consider the diagram
\[ \xymatrix{
\GW^{2i}_0(X)\ar[rr]^-{\simeq} \ar[d]_{\omega(n)^2 \psi^n} && \GW^{2(i+2)}_0(\Tc^{\wedge 2} \wedge X_+) \ar[d]^{\psi^n} \ar[r] & \GW^{2(i+2)}_0(\HP^1 \times \HP^1 \times X) \ar[d]^{\psi^n}\\ 
\GW^{2n(i+2)-4}_0(X) \ar[rr]^-{\simeq} && \GW^{2n(i+2)}_0(\Tc^{\wedge 2} \wedge X_+) \ar[r] &\GW^{2n(i+2)}_0(\HP^1 \times \HP^1 \times X) 
}\]
where the horizontal composites are given by $x \mapsto p^*(x) \cdot (u_1-\tau)(u_2 -\tau)$. Then, for $x\in \GW^{2i}_0(X)$ we have by Proposition~\ref{prop:Adams_in_lambda_ring} and Proposition--Definition~\ref{prop:omega_n}
\[
\psi^n(p^*(x) \cdot (u_1-\tau)(u_2 -\tau)) = p^*(\psi^n(x)) \cdot \psi^n(u_1-\tau)\cdot \psi^n(u_2 -\tau) = \omega(n)^2 \cdot p^*(\psi^n(x)) (u_1-\tau)(u_2 -\tau),
\]
showing that the exterior square in the above diagram commutes. Since the lower right horizontal arrow is injective (e.g.\ by \cite[Lemma~7.6]{Panin10prec}), it follows that the interior left square commutes. By Lemma~\ref{lemm:lift}, this implies that the following diagram commutes
\[ \xymatrix{
\GW^{2i}\ar[rr] \ar[d]_{\omega(n)^2 \psi^n} && \Omega_\Tc^2\GW^{2(i+2)} \ar[d]^{\Omega_\Tc^2\psi^n} \\ 
\GW^{2n(i+2)-4} \ar[rr] && \Omega_\Tc^2\GW^{2n(i+2)}
}\]
which implies the statement by adjunction.
\end{proof}

\begin{cor}
\label{cor:omegai}
Let $i\in\ZZ$ be odd and let $n \in \NN$. Then the diagram
\[
\xymatrix{\Tc^{\wedge 2} \wedge \GW^{2i}\ar[r]\ar[d]_-{\id_{\Tc^{\wedge 2}} \wedge \Psi^n_i} & \GW^{2(i+2)}\ar[d]^-{\Psi^n_{i+2}} \\
\Tc^{\wedge 2} \wedge \GWi{2i}{n}\ar[r] & \GWi{2(i+2)}{n}
}
\]
commutes in $\Ho(\Zh)$, where the upper horizontal arrow is induced by the map \eqref{eq:def_sigma}, and the lower one by \eqref{eq:GWi_bond}.
\end{cor}
\begin{proof}
We have a commutative diagram
\[ \xymatrix{
\Tc^{\wedge 2} \wedge \GW^{2n(i+2)-4} \ar[r] \ar[d] & \Tc^{\wedge 2} \wedge \GWi{2n(i+2)-4}{n} \ar[rrrr]^-{\id_{\Tc^{\wedge 2}} \wedge \Omega_{\Tc}^{\infty}\Sigma_{\Tc}^{i}\omega(n)^{-i-2}} \ar[d] &&&& \Tc^{\wedge 2} \wedge\GWi{2i}{n} \ar[d]\\ 
\GW^{2n(i+2)} \ar[r]& \GWi{2n(i+2)}{n} \ar[rrrr]^-{\Omega_{\Tc}^{\infty}\Sigma_{\Tc}^{i+2}\omega(n)^{-i-2}} &&&& \GWi{2(i+2)}{n}
}\]
Combining this diagram with Proposition~\ref{prop:omegai} yields the corollary, in view of \eqref{eq:def_Psi}.
\end{proof}

\begin{prop}
\label{prop:propuniquelift}
For any $r,n \in \NN$, the natural morphism
\[
\Big[\sOGr^{\wedge r},\sOGr\Big[\frac{1}{n^*}\Big]\Big]_{\SH(\Zh)} \to \lim_{i \text{ odd}} \Big[(\GW^{2i})^{\wedge r},\GWibig{2ir}{n}\Big]_{\Ho(\Zh)}
\]
is bijective.
\end{prop}
\begin{proof}
We use (the proof of) \cite[Theorem 13.1]{Panin10prec} (which applies to $S= \Spec(\Zh)$ by \cite[Theorems 13.2 and 13.3]{Panin10prec}), with the difference that $\BO=\sOGr$ should be replaced by $\sOGr[\frac{1}{n^*}]$, which does not affect any of the arguments appearing in its proof, by \eqref{eq:GW_inverted}. This yields the natural isomorphism
\[
\Big[\sOGr^{\wedge r},\sOGr\Big[\frac{1}{n^*}\Big]\Big]_{\SH(\Zh)} \to \lim_{i \in \NN} \Big[(\GW^{2i})^{\wedge r},\GWibig{2ir}{n}\Big]_{\Ho(\Zh)},
\]
and the proposition follows using a cofinality argument.
\end{proof}

The transition maps in the limit appearing in Proposition~\ref{prop:propuniquelift} are given by the composite
\begin{align*}
\Big[(\GW^{2(i+2)})^{\wedge r},\GWibig{2(i+2)r}{n}\Big]_{\Ho(\Zh)}
&\to  \Big[(\Tc^{\wedge 2} \wedge \GW^{2i})^{\wedge r},\GWibig{2(i+2)r}{n}\Big]_{\Ho(\Zh)} \\ 
&= \Big[(\Sigma_\Tc^2 \Sigma^\infty_\Tc \GW^{2i})^{\wedge r},\Sigma_\Tc^{(i+2)r}\sOGr\Big[\frac{1}{n^*}\Big]\Big]_{\SH(\Zh)}\\
&= \Big[(\Sigma^\infty_\Tc \GW^{2i})^{\wedge r},\Sigma_\Tc^{ir}\sOGr\Big[\frac{1}{n^*}\Big]\Big]_{\SH(\Zh)}\\
&=\Big[(\GW^{2i})^{\wedge r},\GWibig{2ir}{n}\Big]_{\Ho(\Zh)},
\end{align*}
where the first map is given by composition with the map $\Tc^{\wedge 2} \wedge \GW^{2i}\to\GW^{2(i+2)}$ induced by \eqref{eq:GWi_bond}. It thus follows from  Proposition~\ref{cor:omegai} that the family $\Psi_i^n$ of \eqref{eq:def_Psi}, for $i$ odd, defines an element of the limit appearing in Proposition~\ref{prop:propuniquelift} (with $r=1$).

\begin{df}
\label{df:Psi}
For $n \in \NN$, we denote by
\[
\Psi^n \colon \sOGr \to \sOGr\Big[\frac{1}{n^*}\Big].
\]
the morphism of spectra corresponding to the family $\Psi_i^n$ of \eqref{eq:def_Psi}, for $i$ odd, under the bijection of Proposition~\ref{prop:propuniquelift}  (with $r=1$). We call it the \emph{stable $n$-th Adams operation}. 
\end{df}
\begin{rem}
If $X$ is a regular $\Zh$-scheme with structural morphism $p_X \colon X \to \Spec(\Zh))$, we obtain a morphism of spectra 
\[
\Psi^n \colon \sOGr_X = p_X^*\sOGr \to p_X^*\Big( \sOGr\Big[\frac{1}{n^*}\Big]\Big) = (p_X^* \sOGr)\Big[\frac{1}{n^*}\Big] = \sOGr_X\Big[\frac{1}{n^*}\Big].
\]
\end{rem}

For $i\in \ZZ$, let us define
\[
\widetilde{\Psi}_i^n = \Omega^\infty_\Tc \Sigma^i_{\Tc}(\Psi^n) \colon \GW^{2i} \to \GWibig{2ni}n.
\]
Note that, by construction, we have $\Psi_i^n = \widetilde{\Psi}_i^n$ when $i$ is odd. Let us mention that the stable Adams operation has the expected relation to the unstable one, also in even degrees:
\begin{lem}
\label{lemm:Psi_psi}
When $X$ is a smooth $\Zh$-scheme, for any $i\in \ZZ$, the morphism $\GW^{2i}_0(X) \to \GW^{2i}_0(X)[\frac{1}{n^*}]$ induced by $\widetilde{\Psi}_i^n$ equals $\omega(n)^{-i} \psi^n$.
\end{lem}
\begin{proof}
This is true when $i$ is odd, since $\Psi_i^n = \widetilde{\Psi}_i^n$ in this case. Assume that $i$ is even. Let $p \colon X \times \HP^1\to X$ be the projection. We have a commutative diagram
\[ \xymatrix{
\GW^{2i}_0(X)\ar[r]^-\sim \ar[d]_{\widetilde{\Psi}_i^n} & \GW^{2i+2}_0(\Tc \wedge X_+) \ar[d]^{\widetilde{\Psi}_{i+1}^n} \ar[r]& \GW^{2i+2}_0(\HP^1 \times X) \ar[d]^{\widetilde{\Psi}_{i+1}^n}\\ 
\GW^{2i}_0(X)[\frac{1}{n^*}] \ar[r]^-\sim & \GW^{2i+2}_0(\Tc \wedge X_+)[\frac{1}{n^*}] \ar[r] & \GW^{2i+2}_0(\HP^1 \times X) [\frac{1}{n^*}]
}\]
where the horizontal composites are given by $x \mapsto p^*(x) \cdot (u-\tau)$. Now, by the odd case treated above, we have for $x\in \GW^{2i}_0(\HP^1 \times X)$
\begin{align*}
\widetilde{\Psi}^n_{i+1}( p^*(x) \cdot (u-\tau) ) 
&= \omega(n)^{-i-1} \psi^n(p^*(x) \cdot (u-\tau)) && \text{by the odd case}
 \\ 
&= \omega(n)^{-i-1} \cdot p^*\psi^n(x) \cdot \psi^n(u-\tau) && \text{by Proposition~\ref{prop:Adams_in_lambda_ring}}\\
&= p^*(\omega(n)^{-i} \cdot \psi^n(x)) \cdot (u-\tau) && \text{by Proposition--Definition~\ref{prop:omega_n}}.
\end{align*}
The statement then follows from the injectivity of the lower horizontal composite (e.g.\ by \cite[Lemma~7.6]{Panin10prec}).
\end{proof}

\begin{thm}
For any integer $n \in \NN$, the stable Adams operation $\Psi^n \colon \sOGr \to \sOGr[\frac{1}{n^*}]$ is a morphism of ring spectra.
\end{thm}
\begin{proof}
We have first to check that the diagram in $\SH(\Zh)$
\begin{equation}
\begin{gathered}
\xymatrix{
\sOGr\wedge \sOGr\ar[rr]^-{\Psi^n\wedge \Psi^n}\ar[d] && \sOGr[\frac{1}{n^*}]\wedge \sOGr[\frac{1}{n^*}]\ar[d] \\
\sOGr\ar[rr]^-{\Psi^n} && \sOGr[\frac{1}{n^*}]
}
\end{gathered}
\end{equation}
commutes, where the vertical arrows are the multiplications. In view of Proposition~\ref{prop:propuniquelift}, we have to check that the following diagram, in which the vertical maps are induced by the multiplication in the ring spectrum $\sOGr$, commutes in $\Ho(\Zh)$, for any $i\in\ZZ$ odd
\[
\xymatrix@C=8em{
\GW^{2i} \wedge \GW^{2i}\ar[r]^-{\widetilde{\Psi}_i^n\wedge \widetilde{\Psi}_i^n}\ar[d] &  \GWi{2i}{n} \wedge \GWi{2i}{n}\ar[d] \\
\GW^{4i}\ar[r]^{\widetilde{\Psi}_{2i}^n} & \GWi{4i}{n}.}
\]
By Lemma~\ref{lemm:lift} and Lemma~\ref{lemm:Psi_psi} (taking into account \cite[Theorem~11.4]{Panin10prec}), this reduces to the formula, when $X,Y$ are smooth $\Zh$-schemes and $x \in \GW^{2i}_0(X)$, $y \in \GW^{2i}_0(Y)$
\begin{equation}
\label{eq:p_1-p_2}
p_1^*(\omega(n)^{-i} \psi^n(x)) \cdot p_2^*(\omega(n)^{-i} \psi^n(x)) = \omega(n)^{-2i} \psi^n(p_1^*(x) \cdot p_2^*(y)) \in \GW^{4i}_0(X \times Y),
\end{equation}
where $p_1 \colon X \times Y \to X$, $p_2 \colon X \times Y \to Y$ are the projections. But the formula \eqref{eq:p_1-p_2} readily follows from Proposition~\ref{prop:Adams_in_lambda_ring}.

Next, we need to prove the commutativity of the diagram in $\SH(\Zh)$
\[
\xymatrix{\mathbb{S}\ar[r]^-\varepsilon\ar[rd]_-\varepsilon & \sOGr\ar[d]^-{\Psi^n} \\ & \sOGr[\frac{1}{n^*}]}
\]
By adjunction, this reduces to the fact that
\[
\widetilde{\Psi}^n_0(1) =1 \in \GW^0_0(\Spec \Zh),
\]
a consequence of Lemma~\ref{lemm:Psi_psi} and of the fact that $\psi^n(1)=1$.
\end{proof}

\begin{prop}
\label{prop:Psim_Psin}
For any integers $m,n \in \NN$, the composite in $\SH(\Zh)$
\[
\sOGr\xrightarrow{\Psi^n} \sOGr\Big[\frac{1}{n^*}\Big]\xrightarrow{\Psi^m[\frac{1}{n^*}]} \sOGr\Big[\frac{1}{m^*}\Big]\Big[\frac{1}{n^*}\Big] = \sOGr\Big[\frac{1}{(mn)^*}\Big]
\]
is equal to $\Psi^{mn}$. (Here $\Psi^m[\frac{1}{n^*}]$ denotes the image of the morphism $\Psi^m$ under the localisation functor, and the last equality follows from \eqref{eq:m_n_*}.)
\end{prop}
\begin{proof}
For every $i\in \ZZ$, applying the functor $\Omega^\infty_\Tc \Sigma^i_\Tc \colon \SH(\Zh) \to \Ho(\Zh)$ to the morphism $\Psi^m[\frac{1}{n^*}] \colon \sOGr[\frac{1}{n^*}] \to \sOGr[\frac{1}{(mn)^*}]$ yields a morphism
\[
\widetilde{\Psi}^m_i\Big\{\frac{1}{n^*}\Big\} \colon \GWibig{2i}{n} \to \GWibig{2i}{(mn)}.
\]
In view of Proposition~\ref{prop:propuniquelift}, it will suffice to show that, for $i\in \NN$ odd, the composite
\[
\GW^{2i} \xrightarrow{\widetilde{\Psi}^n_i} \GWibig{2i}{n} \xrightarrow{\widetilde{\Psi}^m_i} \GWibig{2i}{(mn)}
\]
equals $\widetilde{\Psi}^{mn}_i$ in $\Ho(\Zh)$. By Lemma~\ref{lemm:lift} and Lemma~\ref{lemm:Psi_psi}, it will then suffice to show that, for each odd $i\in \NN$ and each smooth $\Zh$-scheme $X$, the composite
\[
\GW^{2i}_0(X) \xrightarrow{\omega(n)^{-i}\cdot \psi^n} \GW^{2i}_0(X)\Big[\frac{1}{n^*}\Big] \xrightarrow{\omega(m)^{-i} \cdot \psi^m} \GW^{2i}_0(X)\Big[\frac{1}{(mn)^*}\Big]
\]
equals $\omega(mn)^{-i}\cdot \psi^{mn}$. But this follows from Proposition~\ref{prop:Adams_in_lambda_ring} and Lemma~\ref{lem:omega_psi}.
\end{proof}

\section{Ternary laws for Hermitian \texorpdfstring{$K$}{K}-theory}\label{sec:ternary}
\numberwithin{thm}{section}
\numberwithin{equation}{section}
\renewcommand{\theequation}{\thesection.\alph{equation}}

Recall from \cite[\S 2.3]{Deglise20} that ternary laws are the analogues for Sp-oriented cohomology theories (or spectra) of formal group laws for oriented cohomology theories. In short, the problem is to understand the Borel classes (in the relevant cohomology theory) of the symplectic bundle $U_1\otimes U_2\otimes U_3$ on $\HP^n\times\HP^n\times \HP^n$, where $U_i$ are the universal bundles on the respective factors. The ternary laws permit to compute Borel classes of threefold products of symplectic bundles. At present, there are few computations of such laws, including MW-motivic cohomology and motivic cohomology which are examples of the so-called \emph{additive ternary laws} \cite[Definition 3.3.3]{Deglise20}. In this section, we compute the ternary laws of Hermitian $K$-theory (and thus also of $K$-theory as a corollary), which are not additive.\\

Our first task is to express the Borel classes in Hermitian $K$-theory in terms of the $\lambda$-operations. We will denote by $\sigma_i(X_1,\dots,X_4) \in \ZZ[X_1,\dots,X_4]$ the elementary symmetric polynomials.

\begin{lem}\label{lm:preliminary}
Let $X$ be a $\Zh$-scheme and let $e_1,\ldots,e_4 \in \GW^2_0(X)$ be the classes of rank two symplectic bundles over $X$. Then
\[
\lambda^i(e_1+\dots+e_4)=\begin{cases} \sigma_1(e_1,\ldots,e_4) & \text{ if $i=1$.} \\
\sigma_2(e_1,\ldots,e_4) + 4\gamma & \text{ if $i=2$.} \\
\sigma_3(e_1,\ldots,e_4) + 3\sigma_1(e_1,\ldots,e_4)\gamma & \text{ if $i=3$.} \\
\sigma_4(e_1,\ldots,e_4) + 2\sigma_2(e_1,\ldots,e_4)\gamma + 6\gamma^2 & \text{ if $i=4$.}\end{cases} 
\]
\end{lem}

\begin{proof}
In view of \eqref{eq:lambda_rank_2}, it suffices to expand the product
\[
(1+te_1+\gamma t^2)(1+te_2+\gamma t^2)(1+te_3+\gamma t^2)(1+te_4+\gamma t^2).\qedhere
\]
\end{proof}

\begin{lem}\label{lm:symmetric}
In the ring $\ZZ[x_1,x_2,x_3,x_4,y]$, we have the following equalities:
\[
\sigma_i(x_1-y,\ldots,x_4-y)=\begin{cases} \sigma_1-4y & \text{ if $i=1$,} \\
\sigma_2-3y\sigma_1+6y^2 & \text{ if $i=2$,} \\
\sigma_3-2\sigma_2y+3\sigma_1y^2-4y^3 & \text{ if $i=3$,} \\
\sigma_4-\sigma_3y+\sigma_2y^2-\sigma_1y^3+y^4 & \text{ if $i=4$,}\end{cases}
\]
where $\sigma_i=\sigma_i(x_1,\ldots,x_4)$ for any $i \in \{1,\ldots,4\}$.
\end{lem}

\begin{proof}
Direct computation.
\end{proof}

In the next statement $b_i^{\GW}$ denotes the $i$-th Borel class with values in Hermitian $K$-theory \cite[Definition 8.3]{Panin10pred}.

\begin{prop}\label{prop:borel}
Let $X$ be a $\Zh$-scheme. Let $E$ be a symplectic bundle of rank $8$ on $X$, and $e \in \GW^2_0(X)$ its class. Then we have:
\[
b_i^{\GW}(E)=\begin{cases} e-4\tau & \text{ if $i=1$.} \\
\lambda^2(e)-3\tau e+4(2-3\epsilon)\gamma& \text{ if $i=2$.} \\
\lambda^3(e)-2\tau\lambda^2(e)+3 (1-2\epsilon)\gamma e-8\tau \gamma & \text{ if $i=3$.} \\
\lambda^4(e)- \tau \lambda^3(e)-2\epsilon\gamma\lambda^2(e)-\tau\gamma e+2\gamma^2 & \text{ if $i=4$.} \end{cases}
\]
\end{prop}

\begin{proof}
Using the symplectic splitting principle \cite[\S10]{Panin10pred}, we may assume that $E$ splits as an orthogonal sum of rank two symplectic bundles, whose classes in $\GW^2_0(X)$ we denote by $e_1,\dots,e_4$. The Borel classes $b_i^{\GW}(E)$ are then given by the elementary symmetric polynomials in the elements $e_1-\tau,\ldots,e_4-\tau$, which can be computed using Lemma \ref{lm:symmetric}. For $i=1$, the result is immediate. For $i=2$, we have 
\[
\sigma_2(e_1-\tau,\ldots,e_4-\tau)=\sigma_2(e_1,\ldots,e_4)-3\tau \sigma_1(e_1,\ldots,e_4)+6\tau^2
\]
and $\sigma_2(e_1,\ldots,e_4) = \lambda^2(e) - 4\gamma$ by Lemma \ref{lm:preliminary}. As $\tau^2=2(1-\epsilon)\gamma$, we find 
\[
\sigma_2(e_1-\tau,\ldots,e_4-\tau)= \lambda^2(e) - 4\gamma-3\tau e+12(1-\epsilon)\gamma
\]
proving the case $i=2$. We now pass to the case $i=3$. Using Lemma \ref{lm:symmetric}, we find 
\begin{eqnarray*}
b_3^{\GW}(E) & = & \sigma_3(e_1,\ldots,e_4)-2\tau\sigma_2(e_1,\ldots,e_4)+3\tau^2e-4\tau^3\\
 & = & \lambda^3(e)-3\gamma e-2\tau(\lambda^2(e)-4\gamma)+6(1-\epsilon)\gamma e-16\tau\gamma \\
 & = & \lambda^3(e)+3 (1-2\epsilon)\gamma e-2\tau\lambda^2(e)-8\tau \gamma.
\end{eqnarray*}
In case $i=4$, we have 
\[
b_4^{\GW}(E)=\sigma_4(e_1,\dots,e_4)-\tau \sigma_3(e_1,\dots,e_4)+\tau^2\sigma_2(e_1,\dots,e_4)-\tau^3e+\tau^4.
\]
Using Lemma \ref{lm:preliminary}, we find
\[
\sigma_4(e_1,\dots,e_4) = \lambda^4(e)-2\sigma_2(e_1,\dots,e_4)\gamma-6\gamma^2 = \lambda^4(e)-2 \lambda^2(e)\gamma +2\gamma^2,
\]
\[
\tau\sigma_3(e_1,\dots,e_4)= \tau(\lambda^3(e) -3\gamma e) = \tau \lambda^3(e) -3\tau\gamma e,
\]
\[
\tau^2\sigma_2(e_1,\dots,e_4)=2(1-\epsilon)\gamma\sigma_2(e_1,\dots,e_4)=2(1-\epsilon)\gamma\lambda^2(e) -8(1-\epsilon)\gamma^2.
\]
Since $\tau^3e=4\tau \gamma e$ and $\tau^4=8(1-\epsilon)\gamma^2$, we conclude summing up the previous expressions.
\end{proof}

Our next task is to obtain an explicit formula for the $\lambda$-operations on products of three classes of rank two symplectic bundles, providing a different proof of \cite[Lemma 8.2]{Anan_oper}. It will be useful to have a basis for the symmetric polynomials in three variables $u_1,u_2,u_3$. Following \cite[\S 2.3.3]{Deglise20}, we set, for $i,j,k \in \NN$,
\begin{equation}
\label{eq:sigma_u123}
\sigma(u_1^iu_2^ju_3^k)=\sum_{(a,b,c)}u_1^au_2^bu_3^c
\end{equation}
where the sum runs over the monomials $u_1^au_2^bu_3^c$ in the orbit of $u_1^iu_2^ju_3^k$ under the action of the permutation of the variables $u_1,u_2,u_3$.

\begin{lem}\label{lem:explicit3fold}
Let $X$ be a $\Zh$-scheme, and let $u_1,u_2,u_3 \in \GW^2_0(X)$ be the classes of rank two symplectic bundles on $X$. Then
\[
\lambda^i({u_1}{u_2}{u_3})
=
\begin{cases} 
u_1u_2u_3& \text{if $i=1$.} \\
\sigma(u_1^2u_2^2)\gamma-2\sigma(u_1^2)\gamma^2+4\gamma^3 & \text{if $i=2$.} \\
\sigma(u_1^3u_2u_3)\gamma^2-5u_1u_2u_3\gamma^3 & \text{if $i=3$.} \\
\sigma(u_1^4)\gamma^4+u_1^2u_2^2u_3^2\gamma^3-4\sigma(u_1^2)\gamma^5+6\gamma^6 & \text{if $i=4$.}
\end{cases}
\]
\end{lem} 
\begin{proof}
In view of \eqref{eq:lambda_xyz} and \eqref{eq:lambda_rank_2}, this follows from Lemma \ref{lemm:R_abc}.
\end{proof}

Finally, we are in position to compute the ternary laws of Hermitian K-theory. The computation is obtained by combining Proposition \ref{prop:borel} and Lemma \ref{lem:explicit3fold} (applied to $\gamma^{-1}u_1u_2u_3$). 

\begin{prop}
Let $E_1, E_2, E_3$ be symplectic bundles of rank $2$ on a $\Zh$-scheme $X$. Let $u_1,{u_2},{u_3}$ be their respective classes in $\GW^2_0(X)$. Then the Borel class $b_i^{\GW}(E_1\otimes E_2\otimes E_3) \in \GW^{2i}_0(X)$ equals (using the notation of \eqref{eq:sigma_u123})
\[
\begin{cases}
u_1u_2u_3\gamma^{-1}-4\tau & \text{if $i=1$,} \\
\sigma(u_1^2u_2^2)\gamma^{-1}-2\sigma(u_1^2)-3\tau u_1u_2u_3\gamma^{-1}+12(1-\epsilon)\gamma & \text{if $i=2$,} \\
\sigma(u_1^3u_2u_3)\gamma^{-1}-2(1+3\epsilon)u_1u_2u_3-2\tau\gamma^{-1}\sigma(u_1^2u_2^2)+4\tau\sigma(u_1^2)-16\tau\gamma & \text{if $i=3$,} \\
\sigma(u_1^4)+u_1^2u_2^2u_3^2\gamma^{-1}-4(1-\epsilon)\gamma\sigma(u_1^2)-2\epsilon\sigma(u_1^2u_2^2)-\tau\sigma(u_1^3u_2u_3)\gamma^{-1}+4\tau u_1u_2u_3+8(1-\epsilon)\gamma^2 & \text{if $i=4$.}
\end{cases}
\]
\end{prop}

As a consequence of this proposition, we obtain the explicit expression of the ternary laws  associated to Hermitian $K$-theory (see \cite[Definition 2.3.2]{Deglise20}). We use the notation \eqref{eq:sigma_u123}.
\begin{thm}
The ternary laws $F_i=F_i(v_1,v_2,v_3)$ of Hermitian $K$-theory (over the base $\Spec(\Zh)$) are
\small
\[
F_1=2(1-\epsilon)\sigma(v_1)+ \tau\gamma^{-1} \sigma(v_1v_2)+ \gamma^{-1}v_1v_2v_3, \hspace{10cm} 
\]
\[
F_2=2(1-2\epsilon)\sigma(v_1^2)+2(1-\epsilon)\sigma(v_1v_2)+2\tau\gamma^{-1}\sigma(v_1^2v_2)-3\tau\gamma^{-1} v_1v_2v_3+\gamma^{-1}\sigma(v_1^2v_2^2), \hspace{3cm} 
\]
\[
F_3= 2(1-\epsilon)\sigma(v_1^3)-2(1-\epsilon)\sigma(v_1^2v_2)+8(2-3\epsilon)v_1v_2v_3+\tau \gamma^{-1}\sigma(v_1^3v_2)-2\tau\gamma^{-1} \sigma(v_1^2v_2^2)+3\tau\gamma^{-1}\sigma(v_1^2v_2v_3)+ \gamma^{-1}\sigma(v_1^3v_2v_3),
\]
\[
F_4= \sigma(v_1^4)-2(1-\epsilon)\sigma(v_1^3v_2)+2(1-2\epsilon)\sigma(v_1^2v_2^2)+2(1-\epsilon)\sigma(v_1^2v_2v_3)-\tau\gamma^{-1}\sigma(v_1^3v_2v_3)+2\tau\gamma^{-1}\sigma(v_1^2v_2^2v_3)+\gamma^{-1}\sigma(v_1^2v_2^2v_3^2).
\]
\end{thm}
\begin{proof}
We use the relations $v_i=u_i-\tau$ and the previous theorem. For $b_1$, we find 
\[
u_1u_2u_3=v_1v_2v_3+\tau\sigma(v_1v_2)+\tau^2\sigma(v_1)+\tau^3
\]
and the result follows quite easily from $\tau^2=2(1-\epsilon)\gamma$ and $\tau^3=4\tau\gamma$. For $i=2$, we first compute
\[
\sigma(u_1^2u_2^2)=\sigma(v_1^2v_2^2)+2\tau\sigma(v_1^2v_2)+4(1-\epsilon)\gamma\sigma(v_1^2)+8(1-\epsilon)\gamma\sigma(v_1v_2)+16\tau \gamma\sigma(v_1)+24(1-\epsilon)\gamma^2.
\]
Next,
\[
-2\sigma(u_1^2)=-2\sigma(v_1)^2-4\tau\sigma(v_1)-12(1-\epsilon)\gamma
\]
As $b_2=\sigma(u_1^2u_2^2)\gamma^{-1}-2\sigma(u_1^2)-3\tau u_1u_2u_3\gamma^{-1}+12(1-\epsilon)$, we finally obtain the result for $b_2$.

We now treat the case $i=3$, for which we have 
\[
b_3=\sigma(u_1^3u_2u_3)\gamma^{-1}+(-2-6\epsilon)u_1u_2u_3-2\tau\gamma^{-1}\sigma(u_1^2u_2^2) +4\tau\sigma(u_1^2)-16\tau\gamma
\]
Now,
\[
\sigma(u_1^3u_2u_3)=\sigma(v_1^3v_2v_3)+\tau \sigma(v_1^3v_2)+2(1-\epsilon)\gamma\sigma(v_1^3)+3\tau \sigma(v_1^2v_2v_3)+6(1-\epsilon)\gamma\sigma(v_1^2v_2)+
\]
\[
+12\tau\gamma \sigma(v_1^2)+18(1-\epsilon)\gamma v_1v_2v_3+28\tau\gamma \sigma(v_1v_2)+40(1-\epsilon)\gamma^2 \sigma(v_1) +48 \tau\gamma^2
\]
and we deduce that 
\[
b_3=2(1-\epsilon)\sigma(v_1^3)-2(1-\epsilon)\sigma(v_1^2v_2)+8(2-3\epsilon)v_1v_2v_3+\tau \gamma^{-1}\sigma(v_1^3v_2)-
\]
\[
-2\tau\gamma^{-1} \sigma(v_1^2v_2^2)+3\tau\gamma^{-1}\sigma(v_1^2v_2v_3)+ \gamma^{-1}\sigma(v_1^3v_2v_3).
\]
We conclude with the case $i=4$. The Borel class reads 
\[
b_4=\sigma(u_1^4)+\gamma^{-1}\sigma(u_1^2u_2^2u_3^2)-2\sigma(u_1^2u_2^2)-\tau\gamma^{-1}\sigma(u_1^3u_2u_3)+4\tau u_1u_2u_3 +
\]
\[
+2(1-\epsilon)\sigma(u_1^2u_2^2)-4(1-\epsilon)\gamma\sigma(u_1^2)+8(1-\epsilon)\gamma^2.
\]
First, we note that 
\[
\sigma(u_1^4)=\sigma(v_1^4)+4\tau \sigma(v_1^3)+12(1-\epsilon)\gamma\sigma(v_1^2)+16\tau\gamma\sigma(v_1)+24(1-\epsilon)\gamma^2.
\]
while
\[
u_1^2u_2^2u_3^2=\sigma(v_1^2v_2^2v_3^2)+2\tau \sigma(v_1^2v_2^2v_3)+2(1-\epsilon)\gamma\sigma(v_1^2v_2^2)+8(1-\epsilon)\gamma\sigma(v_1^2v_2v_3)+
\]
\[
+8\tau\gamma\sigma(v_1^2v_2)+32\tau\gamma v_1v_2v_3+8(1-\epsilon)\gamma^2\sigma(v_1^2)+32(1-\epsilon)\gamma^2\sigma(v_1v_2)+32\tau\gamma^2\sigma(v_1)+32(1-\epsilon)\gamma^3.
\]
Using the above, we finally find
\[
b_4=\sigma(v_1^4)-2(1-\epsilon)\sigma(v_1^3v_2)+2(1-2\epsilon)\sigma(v_1^2v_2^2)+2(1-\epsilon)\sigma(v_1^2v_2v_3)-
\]
\[
-\tau\gamma^{-1}\sigma(v_1^3v_2v_3)+2\tau\gamma^{-1}\sigma(v_1^2v_2^2v_3)+\gamma^{-1}\sigma(v_1^2v_2^2v_3^2).
\]
\end{proof}

\begin{rem}
The ternary laws of the spectrum $\mathbf{W}$ representing (Balmer) Witt groups have been computed by Ananyevskiy in \cite[Lemma 8.2]{Anan_oper}. In view of the morphism of ring spectra $\sOGr\to \mathbf{W}$, we may recover this result by setting $1-\epsilon=0$ and $\tau=0$ in the above expression.
\end{rem}

The above theorem yields an expression of the ternary laws of $K$-theory (those can of course be computed more directly). As above, we want to write the Borel classes of threefold products of symplectic bundles in terms of the first Borel classes of the bundles, and we may use the forgetful functor from Hermitian $K$-theory to ordinary $K$-theory. Regarding periodicity, the forgetful functor maps $\tau$ to $2\beta^2$ and $\gamma$ to $\beta^4$, where $\beta$ is the Bott element (of bidegree $(2,1)$). 

\begin{thm}
The ternary laws $F_i=F_i(v_1,v_2,v_3,v_4)$ of $K$-theory are
\[
F_1=4\sigma(v_1)+2\beta^{-2}\sigma(v_1v_2)+\beta^{-4}v_1v_2v_3, \hspace{10cm}
\]
\[
F_2=6\sigma(v_1^2)+4\sigma(v_1v_2)+4\beta^{-2}\sigma(v_1^2v_2)-6\beta^{-2} v_1v_2v_3+\beta^{-4}\sigma(v_1^2v_2^2), \hspace{5cm}
\]
\[
F_3= 4\sigma(v_1^3) -4\sigma(v_1^2v_2)+40v_1v_2v_3+2\beta^{-2}\sigma(v_1^3v_2)-4\beta^{-2}\sigma(v_1^2v_2^2)+6\beta^{-2}\sigma(v_1^2v_2v_3)+\beta^{-4}\sigma(v_1^3v_2v_3),
\]
\[
F_4=\sigma(v_1^4)-4\sigma(v_1^3v_2)+6\sigma(v_1^2v_2^2)+4\sigma(v_1^2v_2v_3)-2\beta^{-2}\sigma(v_1^3v_2v_3)+4\beta^{-2}\sigma(v_1^2v_2^2v_3)+\beta^{-4}v_1^2v_2^2v_3^2. 
\]
\end{thm}

\appendix

\section{\texorpdfstring{$\lambda$}{lambda}-rings}
\label{appendix:lambda}
Here we recall a construction from \cite[V, \S 2.3]{SGA6}; a more accessible exposition can be found in \cite[\S1]{Atiyah-Tall}, where the terminology ``$\lambda$-ring''/``special $\lambda$-ring'' is used instead of ``pre-$\lambda$-ring''/``$\lambda$-ring''. Let $R$ be a commutative ring. One defines a ring $\Lambda(R)$, whose underlying set is $1 + tR[[t]]$. The addition in $\Lambda(R)$ is given by multiplication of power series, while multiplication in  $\Lambda(R)$ is given by the formula
\[
\Big(\sum_{n\in\NN} f_n t^n\Big) \cdot \Big(\sum_{n\in\NN} g_n t^n\Big) = \sum_{n\in\NN} P_n(f_1,\dots,f_n,g_1,\dots,g_n) t^n,
\]
where $P_n$ are certain universal polynomials defined in \eqref{eq:def_P} below. In this ring the neutral element for the addition is the constant power series $1$, and the multiplicative identity is the power series $1+t$. A structure of pre-$\lambda$-ring on $R$ is a morphism of abelian groups
\begin{equation}
\label{def:lambda_t}
\lambda_t = \lambda_t^R \colon R \to \Lambda(R) \quad ; \quad r \mapsto \sum_{n\in \NN} \lambda^n(r) t^n.
\end{equation}
When $R,S$ are pre-$\lambda$-rings, a ring morphism $f\colon R \to S$ is called a morphism of pre-$\lambda$-rings if it commutes with the operations $\lambda^n$, i.e.\ if the following diagram commutes
\[ \xymatrix{
\Lambda(R) \ar[r]^{\Lambda(f)} & \Lambda(S)\\
R\ar[r]^f \ar[u]^{\lambda_t^R} & S \ar[u]_{\lambda_t^S}
}\]

When $R$ is a ring, a pre-$\lambda$-ring structure on $\Lambda(R)$ is defined by setting for $j \in \NN\smallsetminus\{0\}$
\[
\lambda^j\Big(\sum_{n\in \NN} f_n t^n\Big) = \sum_{i\in\NN} Q_{i,j}(f_1,\dots,f_{ij}) t^i,
\]
where $Q_{i,j}$ are certain universal polynomials defined in \eqref{eq:def_Q}. Then $R \mapsto \Lambda(R)$ defines a functor from the category of rings to that of pre-$\lambda$-rings. 

A pre-$\lambda$-ring $R$ is called a $\lambda$-ring if $\lambda_t$ is a morphism of pre-$\lambda$-rings. This amounts to the following relations, for all $n,i,j \in \NN\smallsetminus\{0\}$:
\begin{equation}
\label{eq:L1}
\lambda^n(xy) = P_n(\lambda^1(x),\dots,\lambda^n(x),\lambda^1(y),\dots,\lambda^n(y)) \quad \text{for $x,y \in R$},
\end{equation}
\begin{equation}
\label{eq:L2}
\lambda^i(\lambda^j(z)) = Q_{i,j}(\lambda^1(z),\dots,\lambda^{ij}(z))  \quad \text{for $z \in R$}.
\end{equation}
Note that if $E$ is a subset of $R$ such that \eqref{eq:L1} and \eqref{eq:L2} are satisfied for all $x,y,z \in E$, then \eqref{eq:L1} and \eqref{eq:L2} are satisfied for all $x,y,z$ lying in the subgroup generated by $E$ in $R$.

Note also that if $R$ is a $\lambda$-ring, and $x,y,z \in R$, it follows from Lemma \ref{lemm:R_P} that
\begin{equation}
\label{eq:lambda_xyz}
\lambda^n(xyz) = R_n(\lambda^1(x),\dots,\lambda^n(x),\lambda^1(y),\dots,\lambda^n(y),\lambda^1(z),\dots,\lambda^n(z)),
\end{equation}
where $R_n$ is a polynomial defined in \S\ref{sect:R}.

\begin{lem}
\label{lem:lambda_dim1}
Let $R$ be a commutative ring and $x\in R$. Then in $\Lambda(R)$ we have 
\[
\lambda^1(1+xt)=1+xt \quad \text{and} \quad \lambda^i(1+xt)=0 \text{ for $i>1$}.
\]
\end{lem}
\begin{proof}
This amounts to verifying that $Q_{ij}(x,0,\ldots)=x$ when $i=j=1$, and that $Q_{ij}(x,0,\ldots)=0$ when $i>1$ or $j>1$, which follows at once from \eqref{eq:def_Q} under $U_1 \mapsto x$ and $U_s \mapsto 0$ for $s>0$.
\end{proof}

\begin{lem}
\label{lem:product_dim1}
Let $R$ be a commutative ring and $x\in R$. Let $f_i \in R$ for $i\in \NN$ be such that $f_0=1$. Then
\[
\Big(\sum_{n \in \NN} f_nt^n\Big) \cdot (1+xt)  = \sum_{n \in \NN} f_n x^n t^n \in \Lambda(R).
\]
Moreover, if $x\in R^\times$, then $1+xt$ is invertible in $\Lambda(R)$, and $(1+xt)^i = 1+x^it$ for all $i\in \ZZ$.
\end{lem}
\begin{proof}
The first formula amounts to verifying that $P_n(f_1,\ldots,f_n,x,0,\ldots) = f_n x^n$, which follows from \eqref{eq:def_P} (and \eqref{eq:elementary_sym}) under $V_1 \mapsto x$ and $V_j \mapsto 0$ for $j>1$.
\end{proof}

\begin{lem}
\label{lem:lambdapowerseries}
Let $R$ be a $\lambda$-ring, and consider the ring of Laurent polynomials $R[x^{\pm 1}]$ with coefficients in $R$. Then there exists a unique structure of $\lambda$-ring on $R[x^{\pm 1}]$ such that $R \to R[x^{\pm 1}]$ is a morphism of pre-$\lambda$-rings and $\lambda_t(x)=1+xt$. In addition,
\[
\lambda^n(rx^i) = \lambda^n(r)x^{ni} \quad \text{ for any $r \in R, i\in \ZZ, n \in \NN$}.
\]
\end{lem}
\begin{proof}
Let $S = R[x^{\pm 1}]$. By Lemma \ref{lem:product_dim1}, the element $1+xt \in \Lambda(S)$ is invertible and there exists then a unique pre-$\lambda$-ring structure $\lambda_t \colon S \to \Lambda(S)$ such that $\lambda_t(x)=1+xt$ and $R \to S$ is a morphism of pre-$\lambda$-rings. Consider the diagram
\[\xymatrix{
\Lambda(S)\ar[rrrr]^{\Lambda(\lambda_t^S)}&&&& \Lambda(\Lambda(S))\\
&\Lambda(R)\ar[rr]^{\Lambda(\lambda_t^R)} \ar[lu]&& \Lambda(\Lambda(R)) \ar[ru] &\\ 
&R \ar[rr]^{\lambda_t^R} \ar[u]^{\lambda_t^R} \ar[ld]&& \Lambda(R)\ar[u]_{\lambda_t^{\Lambda(R)}} \ar[rd]& \\
S \ar[rrrr]^{\lambda_t^S} \ar[uuu]^{\lambda_t^S}&&&& \Lambda(S) \ar[uuu]_{\lambda_t^{\Lambda(S)}}
}\]
Using the fact that $\Lambda(R)$ and $\Lambda(S)$ are $\lambda$-rings \cite[Theorem~1.4]{Atiyah-Tall}, we see that all maps are ring morphisms. The interior middle square is commutative because $R$ is a $\lambda$-ring, and the right one because $\Lambda(R) \to \Lambda(S)$ is a morphism of pre-$\lambda$-rings. Commutativity of each of the other three interior squares follows from the fact that $R\to S$ is a morphism of pre-$\lambda$-rings. We conclude that the exterior square is a diagram of $R$-algebras. To verify its commutativity it thus suffices to observe its effect on $x \in S$, which is done using Lemma~\ref{lem:lambda_dim1}. We have proved that $S$ is $\lambda$-ring. The last statement follows from Lemma \ref{lem:product_dim1}.
\end{proof}

\section{Graded rings}\label{sec:gradedrings}

Let $S=S_0\oplus S_1$ be a commutative $\ZZ/2$-graded ring. There is a general procedure to construct a commutative $\ZZ$-graded ring out of $S$, which we now explain. We may consider the ring of Laurent polynomials $S[x^{\pm 1}]$ as a graded ring by setting $\vert x\vert=1$ and $\vert s\vert=0$ for any $s\in S$. We consider the $\ZZ$-graded subgroup $\widehat{S}\subset S[x^{\pm 1}]$ defined by 
\[
\widehat{S}_i:= S_{(i \mod 2)}\cdot x^i, \quad \text{for $i\in \ZZ$}.
\]
It is straightforward to check that $\widehat{S}$ is in fact a $\ZZ$-graded subring of $S[x^{\pm 1}]$, and that the canonical homomorphism of abelian groups $S\to \widehat{S}$ defined by $u \mapsto u x^i$ for $u \in S_i$ and $i=0,1$, has the property that the composite with the projection
\[
S\to \widehat{S}\xrightarrow{\pi}\widehat{S}/(x^2-1)
\]
is an isomorphism of $\ZZ/2$-graded rings.\\

Suppose next that $G$ is an abelian group, and that $S$ is a $G$-graded ring having the structure of a $\lambda$-ring. We will say that $S$ is a $G$-graded $\lambda$-ring if $\lambda^i(r)\in S_{ig}$ for any $i\in \NN$, any $g\in G$ and any $r\in S_g$. As a corollary of Lemma \ref{lem:lambdapowerseries}, we obtain the following result.

\begin{lem}\label{lem:lambdapowerseries2}
Let $S$ be a commutative $\ZZ/2$-graded $\lambda$-ring. Then, the structure of $\lambda$-ring on $S[x^{\pm 1}]$ defined in Lemma \ref{lem:lambdapowerseries} induces a $\lambda$-ring structure on $\widehat{S}$ which turns it into a $\ZZ$-graded $\lambda$-ring. If $r\in \widehat{S}_i$ for some $i\in \ZZ$, there exists a unique $s\in S_{(i \mod 2)}$ such that $r=sx^i$ and we have
\[
\lambda^n(r) = \lambda^n(s)x^{ni} \in \widehat{S}_{ni}.
\]
\end{lem}

\section{Some polynomial identities}\label{sec:polynomial}
\numberwithin{thm}{subsection}
\numberwithin{equation}{subsection}
\renewcommand{\theequation}{\thesubsection.\alph{equation}}

When $U_1,\dots,U_m$ is a series of variables, we denote by $\sigma_n(U) \in \ZZ[U_1,\dots,U_m]$ the elementary symmetric functions, defined by the formula, valid in $\ZZ[U_1,\dots,U_m][t]$,
\begin{equation}
\label{eq:elementary_sym}
\prod_{1 \leq i \leq m} (1+tU_i) = \sum_{n \in \NN} t^n \sigma_n(U).
\end{equation}

\subsection{The polynomials \texorpdfstring{$P_n$}{Pn}}
By the theory of symmetric polynomials, there are polynomials $P_n \in \ZZ[X_1,\dots,X_n,Y_1,\dots,Y_n]$ such that
\begin{equation}
\label{eq:def_P}
\prod_{1 \leq i,j \leq m} (1+tU_i V_j) = \sum_{n \in \NN} t^n P_n(\sigma_1(U) ,\dots, \sigma_n(U),\sigma_1(V),\dots,\sigma_n(V))
\end{equation}
holds in $\ZZ[U_1,\dots,U_m,V_1,\dots,V_m][t]$ for every $m$. 

Let $R$ be a commutative ring. For every $x\in R$, let us define elements $\ell_i(x) \in R$ for each integer $i\geq 1$ by the formula
\begin{equation}
\label{eq:ell}
\ell_i(x) =
\begin{cases}
x & \text{if $i=1$},\\
1 & \text{if $i=2$},\\
0 & \text{if $i >2$}.
\end{cases}
\end{equation}

For elements $a_1,\dots,a_r \in R^{\times}$, we consider the polynomial
\begin{equation}
\label{eq:def_pi}
\pi_{a_1,\dots,a_r}(t) = \prod_{\varepsilon_1,\dots,\varepsilon_r \in \{1,-1\}} (1+ta_1^{\varepsilon_1}\cdots a_n^{\varepsilon_n})\in R[t].
\end{equation}
These polynomials can be expressed inductively as
\begin{equation}
\label{eq:pi_rec}
\pi_{a_1,\dots,a_r}(t) = \pi_{a_1,\dots,a_{r-1}}(ta_r)\cdot \pi_{a_1,\dots,a_{r-1}}(ta_r^{-1}).
\end{equation}

Note that for any $a\in R^\times$
\[
\pi_a(t) = 1 + (a+a^{-1})t +t^2,
\]
and for any $a,b \in R^\times$, setting $x=a+a^{-1}$ and $y =b+b^{-1}$,
\begin{equation}
\label{eq:pi_ab}
\pi_{a,b}(t) = 1 + txy + t^2(x^2+y^2-2) + t^3xy +t^4.
\end{equation}

\begin{lem}
\label{lem:RXY}
Let $R$ be a commutative ring and $x,y \in R$. Then
\[
P_n(\ell_1(x),\dots,\ell_n(x),\ell_1(y),\dots,\ell_n(y)) = \begin{cases}
1 & \text{if $n\in \{0,4\}$},\\
xy & \text{if $n \in \{1,3\}$},\\
x^2+y^2-2 & \text{if $n=2$},\\
0 & \text{if $n>4$}.
\end{cases}
\]
\end{lem}
\begin{proof}
Consider the ring $S=R[a,a^{-1},b,b^{-1}]/(x-a-a^{-1},y-b-b^{-1})$. Then $S$ contains $R$. We have $\sigma_i(a,a^{-1}) = \ell_i(x)$ and $\sigma_i(b,b^{-1}) = \ell_i(y)$ for all $i$, so that, by \eqref{eq:def_P} and \eqref{eq:def_pi}
\[
\pi_{a,b}(t) = \sum_n P_n(\ell_1(x),\dots,\ell_n(x),\ell_1(y),\dots,\ell_n(y)) t^n.
\]
Thus the statement follows from \eqref{eq:pi_ab}.
\end{proof}

\begin{lem}
\label{lem:RB}
Let $R$ be a commutative ring and $n \in \NN \smallsetminus \{0\}$. Then for every $r_1,\dots,r_n \in R$, the element
\[
P_n(r_1,\dots,r_n,\ell_1(B),\dots,\ell_n(B)) - B^nr_n \in R[B]
\]
is a polynomial in $B$ of degree $\leq n-1$.
\end{lem}
\begin{proof}
We may assume that $R=\ZZ[X_1,\dots,X_n]$ and that $r_i=X_i$ for all $i=1,\dots,n$. By algebraic independence of the elementary symmetric polynomials, the ring $R$ is then a subring of $R'=\ZZ[U_1,\dots,U_n]$, via $X_i \mapsto \sigma_i(U)$. The ring $S=R'[B,A,A^{-1}]/(B-A-A^{-1})$ then contains $R'[B]$, and thus also $R[B]$. Since $\sigma_i(A,A^{-1}) = \ell_i(B)$ for all $i$, we have in $S[t]$
\[
\sum_{i=1}^n P_i(\sigma_1(U),\dots,\sigma_n(U),\ell_1(B),\dots,\ell_n(B))t^i  = \prod_{i=1}^n (1+tU_iA)(1+tU_iA^{-1}),
\]
and thus, in $R'[B][t]$,
\[
\sum_{i=1}^n P_i(\sigma_1(U),\dots,\sigma_n(U),\ell_1(B),\dots,\ell_n(B))t^i = \prod_{i=1}^n (1+tU_iB+t^2U_i^2).
\]
Expanding the last product and looking at the $t^n$-coefficients of both sides of the equation, we see that $P_n(\sigma_1(U),\dots,\sigma_n(U),\ell_1(B),\dots,\ell_n(B))$ has leading term $B^n \sigma_n(U)$ as a polynomial in $B$ (in view of \eqref{eq:elementary_sym}).
\end{proof}

\subsection{The polynomials \texorpdfstring{$Q_{i,j}$}{Qij}}\label{subsec:Qij}
By the theory of symmetric polynomials, there are polynomials $Q_{i,j} \in \ZZ[X_1,\dots,X_{ij}]$ (where $i,j \in \NN$) such that
\begin{equation}
\label{eq:def_Q}
\prod_{1 \leq \alpha_1 < \dots < \alpha_j \leq m} (1 + U_{\alpha_1}\cdots U_{\alpha_j} t) = \sum_{i \in \NN} t^i Q_{i,j}(\sigma_1(U),\dots,\sigma_{ij}(U))
\end{equation}
holds in $\ZZ[U_1,\dots,U_m][t]$ for every $m$. For instance, we have $Q_{1,j}=X_j$ for any $j\in \NN \smallsetminus \{0\}$. 

\begin{lem}
\label{lem:RZ}
Let $R$ be a commutative ring and $x\in R$. Then
\[
Q_{i,j}(\ell_1(x),\dots,\ell_{ij}(x)) = 
\begin{cases}
\ell_i(x) & \text{if $j=1$ and $i\neq 0$},\\
1 & \text{if $i=1$ and $j=2$, or if $i=0$},\\
0 & \text{otherwise}. 
\end{cases}
\]
\end{lem}
\begin{proof}
Let $S = R[a,a^{-1}]/(x-a-a^{-1})$. Then $S$ contains $R$. Setting $w_1=a$, $w_2=a^{-1}$ and $w_k =0$ in $S$ for $k>2$, we have $\sigma_k(w) = \ell_k(x)$ for all $k$. Thus for all $j\in \NN$
\[
\sum_{i \in \NN} t^i Q_{i,j}(\ell_1(x),\dots,\ell_{ij}(x)) \overset{\eqref{eq:def_Q}}{=} \prod_{1 \leq \alpha_1 < \dots < \alpha_j \leq m} (1 + w_{\alpha_1}\cdots w_{\alpha_j} t) = \begin{cases}
1+tx + t^2 & \text{if $j=1$},\\
1+t & \text{if $j=2$},\\
1 & \text{otherwise.}\end{cases}\qedhere
\]
\end{proof}

\subsection{The polynomials \texorpdfstring{$R_n$}{Rn}}
\label{sect:R}
By the theory of symmetric polynomials, there are polynomials $R_n \in \ZZ[X_1,\dots,X_n,Y_1,\dots,Y_n,Z_1,\dots,Z_n]$ such that
\[
\prod_{1 \leq i,j,k \leq m} (1+tU_i V_jW_k) = \sum_{n \in \NN} t^n R_n(\sigma_1(U) ,\dots, \sigma_n(U),\sigma_1(V),\dots,\sigma_n(V),\sigma_1(W),\dots,\sigma_n(W))
\]
holds in $\ZZ[U_1,\dots,U_m,V_1,\dots,V_m,W_1,\dots,W_m][t]$ for every $m$.

\begin{lem}
\label{lemm:R_P}
For $n \leq m$, we have in $\ZZ[X_1,\dots,X_m,Y_1,\dots,Y_m,Z_1,\dots,Z_m]$
\[
R_n = P_n(X_1,\dots,X_n,P_1(Y_1,Z_1),\dots,P_n(Y_1,\dots,Y_n,Z_1,\dots,Z_n)).
\]
\end{lem}
\begin{proof}
Observe that, in $\ZZ[U_1,\dots,U_m,V_1,\dots,V_m][t]$,
\[
\prod_{1 \leq i,j \leq m} (1+tU_iV_j) = \prod_{i=1}^m \prod_{j=1}^m(1+tU_iV_j) \overset{\eqref{eq:elementary_sym}}{=} \prod_{i=1}^m \Big(\sum_{n\in\NN} \sigma_n(V) U_i^nt^n\Big).
\]
Since the elements $Y_r = \sigma_r(V)$ for $r=1,\dots,m$ are algebraically independent, in view of \eqref{eq:def_P} it follows that we have in $\ZZ[U_1,\dots,U_m,Y_1,\dots,Y_m][t]$, (writing $Y_s = 0$ for $s>m$)
\begin{equation}
\label{eq:UY}
\sum_{n \in \NN} P_n(\sigma_1(U),\dots,\sigma_n(U),Y_1,\dots,Y_n)t^n = \prod_{i=1}^m \Big(\sum_{n\in \NN} Y_n U_i^nt^n\Big).
\end{equation}
Now in $\ZZ[V_1,\dots,V_m,W_1,\dots,W_m]$, set for any $n\in \NN$,
\[
p_n =  P_n(\sigma_1(V),\dots,\sigma_m(V),\sigma_1(W),\dots,\sigma_m(W)),
\]
so that, in $\ZZ[U_1,\dots,U_m,V_1,\dots,V_m,W_1,\dots,W_m][t]$,
\[
\prod_{1 \leq i,j,k \leq m} (1+tU_iV_jW_k) \overset{\eqref{eq:def_P}}{=} \prod_{i=1}^m \Big(\sum_{n \in \NN} p_n U_i^nt^n\Big) \overset{\eqref{eq:UY}}{=} \sum_{n \in \NN} P_n(\sigma_1(U),\dots,\sigma_n(U),p_1,\dots,p_n)t^n.
\]
Since the elements $X_r = \sigma_r(U), Y_r = \sigma_r(V), Z_r=\sigma_r(W)$ for $r=1,\dots,m$ are algebraically independent, this yields the statement.
\end{proof}

\begin{lem}
\label{lemm:R_abc}
Let $R$ be a commutative ring and $x,y,z \in R$. Then
\begin{align*}
&R_n(\ell_1(x),\dots,\ell_n(x),\ell_1(y),\dots,\ell_n(y),\ell_1(z),\dots,\ell_n(z))  \\
=&\begin{cases}
1 & \text{if $n\in \{0,8\}$},\\
xyz & \text{if $n\in \{1,7\}$},\\
x^2y^2 + x^2z^2+y^2z^2 - 2(x^2+y^2+z^2) +4 & \text{if $n\in \{2,6\}$},\\
x^3yz+xy^3z+xyz^3 - 5xyz & \text{if $n\in \{3,5\}$},\\
x^4 + y^4 + z^4 +x^2y^2z^2-4(x^2+y^2+z^2) +6 & \text{if $n=4$},\\
0 & \text{if $n>8$}.
\end{cases}
\end{align*}
\end{lem}
\begin{proof}
Consider the ring $S=R[a,a^{-1},b,b^{-1},c,c^{-1}]/(x-a-a^{-1},y-b-b^{-1},z-c-c^{-1})$. Then $S$ contains $R$. We have $\sigma_i(a,a^{-1}) = \ell_i(x),\sigma_i(b,b^{-1}) = \ell_i(y), \sigma_i(c,c^{-1}) = \ell_i(z)$ for all $i$. Writing $r_n = R_n(\ell_1(x),\dots,\ell_n(x),\ell_1(y),\dots,\ell_n(y),\ell_1(z),\dots,\ell_n(z))$, we have by definition of $R_n$ and \eqref{eq:def_pi}
\[
\pi_{a,b,c}(t) = \sum_{n \in \NN} r_n t^n \in S[t].
\]
Since $\pi_{a,b,c}(t) = \pi_{a,b}(tc) \cdot \pi_{a,b}(tc^{-1})$ by \eqref{eq:pi_rec}, it follows from \eqref{eq:pi_ab} that $\pi_{a,b,c}(t)$ equals
\[
(1 + txyc + t^2(x^2+y^2-2)c^2 + t^3xyc^3 +t^4c^4)(1 + txyc^{-1} + t^2(x^2+y^2-2)c^{-2} + t^3xyc^{-3} +t^4c^{-4}).
\]
To conclude, we compute the coefficients $r_n$ by expanding the above product. We have $r_0=r_8=1$ and $r_n=0$ for $n>8$, as well as
\[
r_1 =r_7 = xy(c+c^{-1}) = xyz.
\]
Using the fact that $c^2+c^{-2}=z^2-2$, we have
\[
r_2 =r_6 = (x^2+y^2-2)(c^2+c^{-2})+x^2y^2 =   x^2y^2 + x^2z^2+y^2z^2 - 2(x^2+y^2+z^2) +4.
\]
Now $c^3+c^{-3} = z^3 -3z$, hence
\[
r_3 = r_5 =xy(c^3+c^{-3}) + (x^2+y^2-2)xy(c+c^{-1}) = x^3yz+xy^3z+xyz^3 - 5xyz.
\]
Finally $c^4 +c^{-4} = z^4 - 4z^2 +2$, hence
\[
r_4 =  c^4+c^{-4} + x^2y^2(c^2+c^{-2}) + (x^2+y^2-2)^2 = x^4 + y^4 + z^4 +x^2y^2z^2-4(x^2+y^2+z^2) +6.
\]
\end{proof}

\end{document}